\documentclass[11pt]{article}
\usepackage[english]{babel}
\usepackage{amsfonts,amsmath,amsthm,amssymb,graphicx,a4wide}
\usepackage[colorlinks=true]{hyperref}

\theoremstyle{definition}

\newcommand{\dotex}{\frac{d}{dt}}
\newcommand{\ddotex}{\frac{d^2}{dt^2}}

\newcommand{\dv}[2]{{\frac{\partial #1}{\partial #2}}}
\newcommand{\dvv}[2]{{\frac{\partial^2 #1}{\partial {#2}^2}}}

\newcommand{\tr}[1]{\operatorname{Tr}\left(#1\right)}
\newcommand{\trr}[1]{\operatorname{Tr}^2\left(#1\right)}

\newcommand{\EE}[1]{\mathbb{E}\left(#1\right)}

\newcommand{\bra}[1]{\langle #1 |}
\newcommand{\ket}[1]{| #1 \rangle}
\newcommand{\bket}[1]{\langle #1 \rangle}
\newcommand{\proj}[1]{\vert#1\rangle\langle#1\vert}

\newcommand{\NN}{{\mathbb N}}
\newcommand{\RR}{{\mathbb R}}
\newcommand{\CC}{{\mathbb C}}
\newcommand{\HH}{\mathcal{H}}

\newcommand{\ba}{\boldsymbol{a}}
\newcommand{\bA}{\boldsymbol{A}}

\newcommand{\bB}{\boldsymbol{B}}

\newcommand{\bD}{\boldsymbol{D}}
\newcommand{\bH}{\boldsymbol{H}}
\newcommand{\bI}{\boldsymbol{I}}
\newcommand{\bK}{\boldsymbol{K}}
\newcommand{\bL}{\boldsymbol{L}}
\newcommand{\bM}{\boldsymbol{M}}
\newcommand{\bN}{{\boldsymbol{N}}}
\newcommand{\bO}{\boldsymbol{O}}
\newcommand{\bP}{\boldsymbol{P}}

\newcommand{\bU}{\boldsymbol{U}}
\newcommand{\bV}{\boldsymbol{V}}
\newcommand{\bX}{\boldsymbol{X}}

\newcommand{\bSx}{\boldsymbol{\sigma_{\! x}}}
\newcommand{\bSy}{\boldsymbol{\sigma_{\! y}}}
\newcommand{\bSz}{\boldsymbol{\sigma_{\! z}}}
\newcommand{\bSp}{\boldsymbol{\sigma_{\!\text{\bf  +}}}}
\newcommand{\bSm}{\boldsymbol{\sigma_{\!\text{\bf -}}}}

\newcommand{\bNp}{{\bN+\bI}}

\newcommand{\rhoK}{\widetilde{\rho}}
\newcommand{\psiK}{\widetilde{\psi}}
\newcommand{\hK}{\widetilde{h}}
\newcommand{\bMK}{\boldsymbol{\widetilde{M}}}
\newcommand{\bKK}{\boldsymbol{\widetilde{K}}}
\newcommand{\rhoKb}{\widetilde{\rho}_{\infty}}
\newcommand{\rhob}{\overline{\rho}}

\newcommand{\rhoe}{\rho^{\text{\tiny est}}}

\newcommand{\minou}{\text{-}}

\newcommand{\LH}{\mathcal{L}(\HH)}
\newcommand{\PP}{\mathbb{P}}
\newcommand{\Ker}{\operatorname{Ker}}

\newcommand{\etab}{\overline{\eta}}
\newcommand{\thetab}{\overline{\theta}}

\title{Models and  Feedback  Stabilization  of  Open Quantum Systems\footnote{This article  is an extended version of the paper  attached to the  conference given by the  author at  the International Congress of Mathematicians in Seoul, August 13 - 21, 2014.}}
\author{Pierre Rouchon \thanks{Centre Automatique et Syst\`{e}mes, Mines ParisTech, PSL Research University,  60 boulevard Saint-Michel, 75006 Paris. E-mail: pierre.rouchon@mines-paristech.fr }}

\date{January 8, 2014}
\begin{document}
\maketitle
\begin{abstract}
 At the quantum level, feedback-loops have  to  take into account  measurement back-action. We present here  the structure of the  Markovian  models including such back-action and  sketch  two  stabilization  methods: measurement-based feedback  where an open quantum system is stabilized  by a classical  controller; coherent or autonomous  feedback   where  a quantum system is stabilized by a quantum  controller with decoherence  (reservoir engineering).
We begin to explain these models and  methods   for  the photon box experiments realized  in  the group of   Serge  Haroche  (Nobel Prize 2012). We present then  these models and  methods for general  open quantum systems.
\end{abstract}

\paragraph{Classification:}
Primary  93B52, 93D15, 81V10, 81P15;
 Secondary  93C20, 81P68,  35Q84.

\paragraph{Keywords:}
Markov model, open quantum system, quantum filtering,   quantum feedback,  quantum master equation.

\section{Introduction}

\begin{figure}[h]
  \centerline{\includegraphics[width=0.8\textwidth]{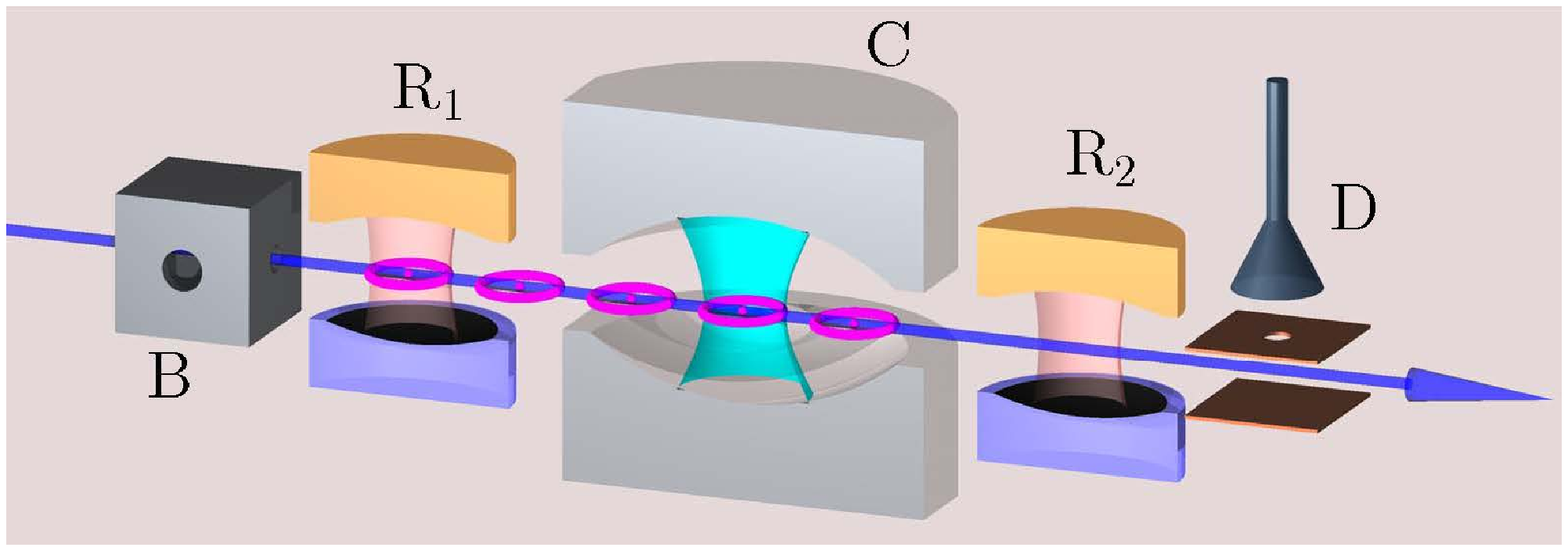}}
  \caption{Scheme of the LBK experiment where photons   are observed via probe atoms.
   The photons in blue  are trapped between the two mirrors of the cavity $C$. They are  probed by  two-level  atoms (the small pink torus) flying out the preparation box $B$, passing  through the cavity  $C$ and measured in $D$. Each atom is manipulated  before and after $C$ in  Ramsey cavities $R_1$ and  $R_2$, respectively.   It is  finally  detected in $D$  either in ground state $\ket g$   or in excited state   $\ket e$.    }\label{fig:LKBsetup0}
\end{figure}

Serge Haroche has obtained  the Physics Nobel Prize in 2012 for a series of  crucial experiments on observations and manipulations of photons with atoms.  The  book~\cite{haroche-raimondBook06}, written with Jean-Michel Raimond,  describes  the physics (Cavity Quantum Electro-Dynamics, CQED)  underlying these  experiments  done  at  Laboratoire Kastler Brossel (LKB).  These  experimental setups,  illustrated on figure~\ref{fig:LKBsetup0} and named in the sequel "the LKB photon box",   rely on  fundamental  examples of open quantum systems constructed with harmonic oscillators and qubits. Their  time evolutions   are captured by stochastic dynamical models  based  on three features,  specific to the quantum world and listed below.
\begin{enumerate}

  \item  The state of a quantum system is described either by the  wave function  $\ket{\psi}$ a vector of length one   belonging to some separable Hilbert space $\HH$ of finite or infinite dimension, or,  more generally, by the density operator $\rho$ that is a non-negative  Hermitian operator on $\HH$ with trace one. When the system can be described by a wave function $\ket{\psi}$ (pure state), the density operator  $\rho$ coincides with  the orthogonal projector on  the line spanned by $\ket{\psi}$ and   $\rho =\ket{\psi}\bra{\psi}$ with usual Dirac notations. In general the rank of $\rho$ exceeds one, the state is then mixed and cannot be described by a wave function.  When the system is closed, the  time evolution of $\ket{\psi}$ is  governed by the Schr\"{o}dinger equation
      \begin{equation}\label{eq:Schrodinger}
        \dotex \ket\psi = - \tfrac{i}{\hbar} \bH \ket\psi
      \end{equation}
      where $\bH$ is the system Hamiltonian, an Hermitian operator on $\HH$  that could possibly depend on time $t$ via some time-varying parameters (classical control inputs).
       When the system is closed, the  evolution of $\rho$ is  governed by the Liouville/von-Neumann equation
      \begin{equation}\label{eq:Liouville}
        \dotex \rho = -\tfrac{i}{\hbar}\big[\bH,\rho\big] =  -\tfrac{i}{\hbar} \big(\bH\rho - \rho\bH\big)
        .
      \end{equation}

  \item Dissipation and irreversibility has its origin in the "collapse of the wave packet"  induced by the measurement. A measurement on the quantum system of state $\ket{\psi}$ or $\rho$  is associated  of  an observable  $\bO$, an Hermitian operator on $\HH$,  with spectral decomposition $\sum_{\mu} \lambda_\mu \bP_\mu$:  $\bP_\mu$ is the orthogonal projector on the eigen-space associated to the eigen-value $\lambda_\mu$.  The measurement process attached to $\bO$ is assumed to be instantaneous and obeys to the following rules:
        \begin{itemize}
          \item the measurement outcome $\mu$ is obtained  with probability $\PP_\mu=\bra{\psi} \bP_\mu\ket\psi$ or $\PP_\mu=\tr{\rho \bP_\mu}$,  depending on the state $\ket{\psi}$ or $\rho$ just before the measurement;
           \item just after the measurement process, the quantum state is changed to $\ket{\psi}_+$ or $\rho_+$ according to the mappings
               $$
  \ket\psi\mapsto\ket{\psi}_+ = \frac{\bP_\mu \ket{\psi}}{\sqrt{\bra\psi \bP_\mu\ket\psi}}\quad  \text{ or }\quad
  \rho\mapsto \rho_+ = \frac{\bP_\mu \rho \bP_\mu}{\tr{\rho \bP_\mu}}
  $$
  where $\mu$  is  the observed  measurement outcome. These mappings  describe the measurement back-action and have  no classical counterpart.
        \end{itemize}

  \item Most systems are composite systems built with several sub-systems.  The quantum states of such composite systems live in the tensor product of the Hilbert spaces of each sub-system. This is a crucial difference with  classical composite systems where the state space is built  with Cartesian products. Such tensor products have important implications such   as  entanglement with   existence of non separable states. Consider a  bi-partite   system made of two sub-systems: the sub-system of interest $S$  with Hilbert space $\HH_S$  and the measured sub-system $M$ with Hilbert space $\HH_M$.  The quantum state of this bi-partite system $(S,M)$ lives in $\HH=\HH_S\otimes \HH_M$. Its Hamiltonian $\bH$ is constructed with the Hamiltonians of the  sub-systems, $\bH_S$ and $\bH_M$, and an interaction  Hamiltonian $ \bH_{int}$ made of a sum of tensor products of operators (not necessarily Hermitian)  on $S$ and $M$:
      $$
      \bH=\bH_S\otimes\bI_M+\bH_{int}+\bI_S\otimes \bH_M
      $$
      with $\bI_S$ and $\bI_M$  identity operators on $\HH_S$ and $\HH_M$, respectively.
      The measurement operator $\bO=\bI_S \otimes \boldsymbol{O}_M $ is here  a simple tensor product of  identity on $S$ and the Hermitian operator $\bO_M$ on $\HH_M$, since only $M$ is directly measured. Its spectrum is  degenerate: the multiplicities of the  eigenvalues are  necessarily greater or equal to the dimension of $\HH_S$.

\end{enumerate}

This paper  shows that, despite  different mathematical  formulations,  dynamical models describing open quantum systems admit the same structure, essentially  given by the Markov model~\eqref{eq:GenMarkovChain}, and  directly  derived  from  the three quantum features listed here above.  Section~\ref{sec:photonBox} explains the construction of such  Markov models for  the LKB photon box and its stabilization  by measurement-based  and coherent feedbacks. These stabilizing feedbacks rely on control Lyapunov functions, quantum filtering and reservoir engineering.  The  next  sections  explain these models and  methods for general  open quantum systems. In section~\ref{sec:DiscreteTime} (resp. section~\ref{sec:ContinuousTime})  general discrete-time  (resp. continuous-time) systems are considered. In appendix, operators,   key states and formulae are  presented    for  the  quantum  harmonic oscillator and for the qubit, two important  quantum systems.  These notations are  used and not explicitly  recalled  throughout sections~\ref{sec:photonBox},  \ref{sec:DiscreteTime} and~\ref{sec:ContinuousTime}.

\section{The LKB photon box} \label{sec:photonBox}

\subsection{The ideal  Markov  model}

The LKB photon box of figure~\ref{fig:LKBsetup0}, a bi-partite system with the photons as first  sub-system  and the probe atom as second sub-system,   illustrates in an almost perfect and  fundamental way   the three quantum features listed in the introduction section.  This system is a discrete time system with sampling period $\tau$ around $80~\mu s$, the time interval between probe  atoms.  Step $k\in\NN$ corresponds to time $t=k\tau$. At $t=k\tau$,  the photons are assumed to be  described by the wave function $\ket{\psi}_k$ of an harmonic oscillator (see appendix~\ref{ap:oscillator}).  At $t=k\tau$, the probe atom number $k$, modeled as a qubit (see appendix~\ref{ap:qubit}),  gets  outside the box $B$   in ground state $\ket{g}$.   Between $t\in[k\tau,(k+1)\tau[$, the wave function $\ket\Psi$ of this composite system,  photons/atom number $k$, is governed by a Schr\"{o}dinger evolution
$$
\dotex \ket{\Psi} = - \frac{i}{\hbar} \bH \ket\Psi
$$
with starting condition $\ket\Psi_{k\tau} = \ket{\psi}_k \otimes \ket{g}$ and
where $\bH$ is the photons/atom  Hamiltonian depending possibly on $t$.   Appendix~\ref{ap:JC} presents   typical  Hamiltonians in the resonant and  dispersive cases.  We have  thus a propagator between $t=k\tau$ and  $t=(k+1)\tau^-$, $U_{(k\tau,(k+1)\tau^-)}$,  from which we get $\ket\Psi$ at time $t=(k+1)\tau^-$, just before  detector $D$ where the energy of the atom is measured via $\bO=\bI_S\otimes\bSz$.
The following relation,
$$
\ket\Psi_{(k+1)\tau^-} = U_{(k\tau,(k+1)\tau^-)} \ket\psi_k \otimes \ket g
\triangleq \bM_g \ket\psi_k \otimes \ket g +  \bM_e \ket\psi_k \otimes \ket e
,
$$
valid for any $\ket\psi_k$,  defines the measurement operators $\bM_g$ and $\bM_e$ on  the Hilbert space of the photons $\HH_S$.  Since, for all $\ket\psi_k$, $\ket\Psi_{(k+1)\tau^-}$ is of length $1$, we have necessarily $\bM_g^\dag \bM_g + \bM_e^\dag \bM_e=\bI_S$.
At time $t=(k+1)\tau^-$,  we measure $\bO = \lambda_e\bI_S\otimes \proj{e} + \lambda_g\bI_S\otimes\proj{g}$ with two highly degenerate eigenvalues $\lambda_e=1$, $\lambda_g=-1$ of eigenspaces $\HH_S \otimes \ket{e}$ and $\HH_S \otimes \ket{g}$, respectively. According to the measurement quantum rules, we can get only two outcomes $\mu$, either $\mu=g$ or $\mu=e$.  With outcome  $\mu$, just after the measurement, at time $(k+1)\tau$  the quantum state $\ket\Psi$ is changed to
$$
\ket\Psi_{(k+1)\tau^-} \mapsto \ket\Psi_{(k+1)\tau}= \frac{\bM_\mu \ket\psi_k }{\sqrt{\bra{\psi_k}\bM_\mu^\dag \bM_\mu\ket{\psi_k}}} \otimes \ket{\mu}.
$$
Moreover the probability to get $\mu$ is $\bra{\psi_k}\bM_\mu^\dag \bM_\mu\ket{\psi_k}$. Since $\ket\Psi_{(k+1)\tau}$ is now a simple tensor product (separate state), we can forget the atom number $k$ and summarize the evolution of the photon wave function between $t=k\tau$ and $t=(k+1)\tau$ by the following Markov process
$$
     \ket\psi_{k+1} = \left\{
                   \begin{array}{ll}
                     \frac{\bM_g \ket\psi_k}{\sqrt{\bra{\psi_k}\bM_g^\dag\bM_g\ket{\psi_k}}},
              & \hbox{with probability $\bra{\psi_k}\bM_g^\dag\bM_g\ket{\psi_k}$;} \\
                  \frac{\bM_e \ket\psi_k}{\sqrt{\bra{\psi_k}\bM_e^\dag\bM_g\ket{\psi_k}}},
              & \hbox{with probability $\bra{\psi_k}\bM_e^\dag\bM_e\ket{\psi_k}$.}
                   \end{array}
                 \right.
$$
More generally, for an arbitrary quantum state $\rho_k$ of the photons at step $k$,  we have
\begin{equation}\label{eq:MarkovChain}
      \rho_{k+1} = \left\{
                   \begin{array}{ll}
                     \frac{\bM_g \rho_k \bM_g^\dag}{\tr{\bM_g \rho_k\bM_g^\dag}},
              & \hbox{with probability $p_g(\rho_k) =\tr{\bM_g \rho_k\bM_g^\dag}$;} \\
                  \frac{\bM_e \rho_k\bM_e^\dag}{\tr{\bM_e \rho_k\bM_e^\dag}},
              & \hbox{with probability $p_e(\rho_k) =\tr{\bM_e \rho_k \bM_e^\dag}$.}
                   \end{array}
                 \right.
\end{equation}
The  measurement operators $\bM_g$ and $\bM_e$ are  implicitly defined by the Schr\"{o}dinger propagator between $k\tau$ and $(k+1)\tau$. They always satisfy $\bM_g^\dag \bM_g + \bM_e^\dag \bM_e=\bI_S$.

\subsection{Quantum Non Demolition (QND) measurement}

For a well tuned composite evolution $U_{(k\tau,(k+1)\tau^-)}$ (see~\cite{haroche-raimondBook06}) with a dispersive interaction, one get the following measurement operators, functions of the photon-number operator~$\bN$,
\begin{equation}\label{eq:MgMeDisp}
  \bM_g= \cos\left(\tfrac{\phi_0 \bN + \phi_R }{2}\right), \quad \bM_e= \sin\left(\tfrac{\phi_0 \bN + \phi_R }{2}\right)
\end{equation}
where $\phi_0$ and $\phi_R$ are  tunable real parameters. The Markov process~\eqref{eq:MarkovChain} admits then  a lot of interesting  properties  characterizing  QND measurement.
\begin{itemize}
  \item For  any function $g:\RR\mapsto \RR$, $V_g(\rho)=\tr{g(\bN) \rho}$ is a martingale:
$$
\EE{V_g(\rho_{k+1})~/~\rho_k} = V_g(\rho_k)
$$
where  $\EE{x~/~y}$ stands for conditional expectation of $x$ knowing $y$.
This results from elementary properties of the trace and from  the commutation of  $\bM_g$ and  $\bM_e$ with  $\bN$.

\item For any integer $\bar n$, the  photon-number state $\ket{\bar n}\bra{\bar n}$ ($\bar n\in\NN$) is a steady-state: any realization  of~\eqref{eq:MarkovChain} starting from $\rho_0=\ket{\bar n}\bra{\bar n}$ is constant:  $\forall k\geq 0$, $\rho_k \equiv \ket{\bar n}\bra{\bar n}$.

\item When $(\phi_R,\phi_0,\pi)$ are $\mathbb{Q}$-independent, there is no other steady state than these photon-number states.   Moreover, for any initial density operator $\rho_0$ with a finite photon-number  support ($\rho_0 \ket{m}=0$ for $m$ large enough), the probability that $\rho_k$ converges towards the steady state $\ket{\bar n}\bra{\bar n}$ is $\tr{\ket{\bar n}\bra{\bar n}\rho_0}= \bra{\bar n}\rho_0\ket{\bar n}$.  Since $\tr{\rho_0}=1= \sum_{\bar n\in\NN} \bra{\bar n}\rho_0\ket{\bar n}$, the Markov process~\eqref{eq:MarkovChain} converges almost surely towards a photon-number state, whatever its initial state $\rho_0$ is.
\end{itemize}
The proof of this convergence result is essentially based on a Lyapunov function, a super-martingale,
$V(\rho)=-\sum_{n\in\NN} \bket{n\left|\rho \right|n}^2$. Simple computations yield
$$
\EE{ V(\rho_{k+1})~/~\rho_k} = V(\rho_k) - Q(\rho_k)
$$
where $Q(\rho) \geq 0$ is given by the following formula
\begin{multline*}
  Q(\rho)=
  \tfrac{\left(\sum_{n'}\cos^2\left(\tfrac{\phi_0 n'+ \phi_R }{2}\right) \bra{n'}\rho\ket{n'}\right)\left(\sum_{n'}\sin^2\left(\tfrac{\phi_0 n'+ \phi_R }{2}\right) \bra{n'}\rho\ket{n'}\right) }{4}
  \\
 \left( \sum_{n\in\NN}
  \left(
  \tfrac{\cos^2\left(\tfrac{\phi_0 n+ \phi_R }{2}\right) \bra{n}\rho\ket{n}}{\sum_{n'} \cos^2\left(\tfrac{\phi_0 n'+ \phi_R }{2}\right) \bra{n'}\rho\ket{n'}}
  -
 \tfrac{\sin^2\left(\tfrac{\phi_0 n+ \phi_R }{2}\right) \bra{n}\rho\ket{n}}{\sum_{n'} \sin^2\left(\tfrac{\phi_0 n'+ \phi_R }{2}\right) \bra{n'}\rho\ket{n'}}
  \right)^2\right)
.
\end{multline*}
Since  $(\phi_0,\phi_R,\pi)$ are $\mathbb{Q}$-independent, $Q(\rho)=0$ implies that,  for some $\bar n\in\NN$, $\rho=\ket{\bar n}\bra{\bar n}$.  One concludes then   with usual probability and compactness   arguments~\cite{kushner-71}, despite the fact that the underlying Hilbert space is of infinite dimension. Other and also  more precise results can be found in~\cite{BauerIHP2013}.

\subsection{Stabilization of photon-number states by feedback}
Take $\bar n\in\NN$. With measurement operators~\eqref{eq:MgMeDisp}, the Markov process~\eqref{eq:MarkovChain} admits  $\bar\rho=\ket{\bar n}\bra{\bar n}$ as steady state. We describe here the  measurement-based feedback  (quantum-state  feedback)  implemented experimentally in~\cite{sayrin-et-al:nature2011}  and that     stabilizes $\bar\rho$. Here the scalar classical  control input  $u$ consists in applying,  just after the  atom measurement in $D$,  a coherent displacement of tunable amplitude $u$. This yields the following control Markov process
\begin{equation}\label{eq:MarkovChainControl}
\rho_{k+1}=
\left\{
  \begin{aligned}
  &\frac{\bD_{u_k}\bM_g\rho_k \bM_g^\dag \bD_{u_k}^\dag}{\tr{\bM_g\rho_k \bM_g^\dag}}\quad \text{$y_k=g$  with probability } p_{g,k}=\tr{\bM_g\rho_k \bM_g^\dag}&\\
  &\frac{\bD_{u_k}\bM_e\rho_k \bM_e^\dag \bD_{u_k}^\dag}{\tr{\bM_e\rho_k \bM_e^\dag}}\quad \text{$y_k=e$  with probability } p_{e,k}=\tr{\bM_e\rho_k \bM_e^\dag}&
    \end{aligned}
  \right.
\end{equation}
where $u_k\in\RR$ is the control at step $k$, $\bD_{u}=e^{u \ba^\dag - u\ba}$ is the displacement of amplitude $u$ (see appendix~\ref{ap:oscillator}) and $y_k$ is the measurement outcome  at step  $k$.

The stabilization  of $\bar\rho$ is based on a state-feedback  function $f$, $u= f(\rho)$,  such that almost all  closed-loop trajectories of~\eqref{eq:MarkovChainControl} with $u_k=f(\rho_k)$ converge towards $\bar\rho$ for any initial condition $\rho_0$.  The construction of $f$ exploits the open-loop martingales $\tr{ g(\bN)\rho}$ to construct the following strict control Lyapunov function:
$$
V_\epsilon(\rho)= \sum_{n} \left( -\epsilon \bket{n\left|\rho \right|n}^2 +
\sigma_n\bra{n}\rho\ket{n}\right)
$$
where $\epsilon  >0$ is small enough and
$$
\sigma_n = \left\{
             \begin{array}{ll}
              \tfrac{1}{4} +  \sum_{\nu=1}^{\bar n} \tfrac{1}{\nu}-\tfrac{1}{\nu^2}, & \hbox{if } n=0; \\
               \sum_{\nu=n+1}^{\bar n} \tfrac{1}{\nu}-\tfrac{1}{\nu^2}, & \hbox{if } n\in[1,\bar n-1]; \\
               0, & \hbox{if } n=\bar n; \\
               \sum_{\nu=\bar n+1}^{n} \tfrac{1}{\nu}+\tfrac{1}{\nu^2}, & \hbox{if } n \in[\bar n +1, +\infty].
             \end{array}
           \right.
$$
The weight $\sigma_n$ are all non negative, $n\mapsto \sigma_n$ is strictly decreasing (resp. increasing)  for $n\leq \bar n$ (resp. $n\geq \bar n$) and minimum for $n=\bar n$.
The feedback law $u=f(\rho)$  is obtained by choosing $u$ such that the   expectation value of $V_{\epsilon}(\rho_{k+1})$, knowing $\rho_k=\rho$ and $u_k=u$,  is as small as possible:
$$
u=f(\rho)=:\underset{\upsilon\in[-\bar u,\bar u]}{\text{Argmin}}\quad  V_\epsilon\Big(\bD_\upsilon\left( \bM_g\rho \bM_g^\dag + \bM_e\rho \bM_e^\dag \right) \bD_\upsilon^\dag\Big)
$$
where $\bar u >0$ is some prescribed  bound on $|u|$. Such a feedback law achieves global stabilization since, in closed-loop,  the Lyapunov function is  strict:
$$
\forall \rho\neq \ket{\bar n}\bra{\bar n}, \quad  V_\epsilon\Big(\bD_{f(\rho)}\left( \bM_g\rho \bM_g^\dag + \bM_e\rho \bM_e^\dag \right)\bD_{f(\rho)}^\dag\Big) <  V_\epsilon\Big(\rho\Big)
.
$$
Formal convergence proofs  can be found in~\cite{AminiSDSMR2013A} for any finite dimensional approximations resulting from a truncation to a finite number of photons and in~\cite{SOMARAJU2013} for  the infinite dimension.

\subsection{A more realistic Markov  model with detection  errors}

The experimental implementation of the above feedback law \cite{sayrin-et-al:nature2011} has to cope with several sources  of imperfections. We focus here on measurement errors and show how the Markov process has to be changed to take into account these errors. Assume that we know the  detection error rates characterized by $\mathbb{P}(y=e/\mu=g)=\eta_g\in[0,1]$ (resp. $\mathbb{P}(y=g/\mu=e)=\eta_e\in[0,1]$)  the probability of erroneous  assignation to $e$ (resp. $g$) when the atom collapses in $g$ (resp. $e$). Without error, the quantum state $\rho_k$ obeys to~\eqref{eq:MarkovChain}.   A direct application of  Bayes law provides the expectation of $\rho_{k+1}$, knowing $\rho_k$ and the  effective  detector signal  $y_k$, possibly corrupted by a detection error.
When $y_k=g$, this expectation value is given by
$
\frac{ (1-\eta_g)\bM_g\rho_k \bM_g^\dag + \eta_e \bM_e\rho_k \bM_e^\dag }{\tr{ (1-\eta_g)\bM_g\rho_k \bM_g^\dag + \eta_e \bM_e\rho_k \bM_e^\dag}}
$
and, when $y_k=e$,  by
$
\frac{\eta_g \bM_g\rho_k \bM_g^\dag +  (1-\eta_e)\bM_e\rho_k \bM_e^\dag  }{\tr{\eta_g \bM_g\rho_k \bM_g^\dag +  (1-\eta_e)\bM_e\rho_k \bM_e^\dag}}
.
$
Moreover the probability to get  $y_k=g$ is {\small $\tr{ (1-\eta_g)\bM_g\rho_k \bM_g^\dag + \eta_e \bM_e\rho_k \bM_e^\dag}$} and to get $y_k=e$ is
 {\small $ \tr{\eta_g \bM_g\rho \bM_g^\dag +  (1-\eta_e)\bM_e\rho \bM_e^\dag}$}. This means that
the  Markov process~\eqref{eq:MarkovChain} must be changed to
\begin{equation}\label{eq:MarkovChainError}
      \rho_{k+1} = \left\{
                   \begin{array}{ll}
                     \frac{ (1-\eta_g)\bM_g\rho_k \bM_g^\dag + \eta_e \bM_e\rho_k \bM_e^\dag }{\tr{ (1-\eta_g)\bM_g\rho_k \bM_g^\dag + \eta_e \bM_e\rho_k \bM_e^\dag}}
              & \hbox{ when $y_k=g$,} \\[1.em]
                     \frac{\eta_g \bM_g\rho_k \bM_g^\dag +  (1-\eta_e)\bM_e\rho_k \bM_e^\dag  }{\tr{\eta_g \bM_g\rho_k \bM_g^\dag +  (1-\eta_e)\bM_e\rho_k \bM_e^\dag}}
              & \hbox{ when $y_k=e$,}
                   \end{array}
                 \right.
\end{equation}
with
{\small $
 \tr{ (1-\eta_g)\bM_g\rho_k \bM_g^\dag + \eta_e \bM_e\rho_k \bM_e^\dag}$} and {\small $\tr{\eta_g \bM_g\rho_k \bM_g^\dag +  (1-\eta_e)\bM_e\rho_k \bM_e^\dag}
$}
being the  probabilities  to detect $y_k=g$ and $e$, respectively.
The quantum state $\rho_k$ is thus a conditional state: it is the expectation value of the projector associated to  the photon  wave function at step $k$, knowing   its value at step $k=0$ and the detection outcomes $(y_0,\ldots,y_{k-1})$.

All other experimental imperfections including decoherence can be treated in the same way (see, e.g., \cite{dotsenko-et-al:PRA09,somaraju-et-al:acc2012}) and yield to a quantum state governed by a Markov  process with a similar structure. In fact all usual models of open quantum systems admit the same structure, either in discrete-time (see section~\ref{sec:DiscreteTime}) or  in continuous-time (see section~\ref{sec:ContinuousTime}).

\subsection{The real-time stabilization algorithm}
Let us give more details on  the real-time implementation used in~\cite{sayrin-et-al:nature2011} of this quantum-state  feedback.
The sampling period $\tau$ is around $80~\mu s$.  The controller set-point is an integer  $\bar n$ labelling  the steady-state $\bar \rho=\ket{\bar n}\bra{\bar n}$  to be stabilized.  At  time step $k$, the real-time computer
   \begin{enumerate}
      \item  reads $y_k$ the measurement outcome  for probe atom $k$;
      \item  updates the quantum state  from previous step value $\rho_{k-1}$ to $\rho_{k}$ using  $y_k$ and  a Markov model slightly more complicated but of same structure as~\eqref{eq:MarkovChainError}; this update corresponds to  a quantum filter (see subsection~\ref{ssec:filtering}).
      \item computes $u_k$ as  $f(\rho_{k})$ (state feedback) where $f$ results from  minimizing the expectation of the control Lyapunov function $V_{\epsilon}(\rho)$ at step $k+1$, knowing $\rho_{k}$;
       \item send via  an  antenna a  micro-wave pulse  calibrated to obtain the displacement  $\bD_{u_{k}}$ on the photons.
   \end{enumerate}
All the details of this quantum feedback  are given in~\cite{sayrin:thesis}. In particular, the Markov  model takes into account several  experimental imperfections such as  finite life-time of the photons (around $1/10~s$) and a delay of $5$~steps  in the feedback loop. Convergence results related to this feedback scheme  are given in~\cite{AminiSDSMR2013A}.

\subsection{Reservoir engineering stabilization of Schr\"{o}dinger cats} \label{ssec:LKBcat}

It is possible to stabilize the  photons trapped in cavity $C$ (figure~\ref{fig:LKBsetup0}) without any such measurement-based  feedback, just by well tuned interactions with the  probe atoms and without measuring  them in $D$. Such kind of stabilization, known as reservoir engineering~\cite{PoyatCZ1996PRL}, can be seen as a generalization of  optical pumping techniques~\cite{Kastl1967S}.  Such stabilization methods are illustrative of coherent (or autonomous) feedback where the controller is an open quantum system.  In~\cite{sarlette-et-al:PRL2011}, a realistic implementation of such passive stabilization method is proposed.   It  stabilizes a coherent superposition of classical  photon-states with opposite phases,  a Schr\"{o}dinger phase-cats with wave functions  of the form $(\ket{\alpha} + i \ket{\minou\alpha})/\sqrt{2}$, where $\ket{\alpha}$ is  the coherent state of amplitude $\alpha\in\RR$.  We explain here the convergence analysis of such passive stabilization using the notations and operator definitions  given in appendix~\ref{ap:oscillator}.

The atom entering the cavity~$C$  is prepared through $R_1$ in a partially excited state $\cos(u/2) \ket{g} + \sin(u/2) \ket e$ with $u\in[0,\pi/2[$ (south hemisphere of the Bloch sphere). Its interaction with the photons is  first dispersive with positive detuning during its entrance,  then resonant in the cavity middle   and finally dispersive with negative detuning when leaving the cavity. The resulting measurement operators $\bM_g$ and $\bM_e$ appearing in~\eqref{eq:MarkovChain} admit then the following form (see \cite{sarlette-et-al:PRA2012} for detailed derivations):
$$
\bM_g=  e^{- i \hK(\bN)} \bMK_g e^{ i \hK(\bN)}, \quad \bM_e=  e^{- i \hK(\bN)} \bMK_e e^{ i \hK(\bN)}
$$
with $n\mapsto \hK(n)$  a real   function, with $\bI$ standing for $\bI_S$, with
\begin{align*}
  \bMK_g &=  \cos(\tfrac{u}{2}) \cos\left(\tfrac{\theta(\bN)}{2}\right) + \epsilon\sin(\tfrac{u}{2}) \frac{\sin\left(\tfrac{\theta(\bN)}{2}\right)}{\sqrt{\bN}} \, \ba^\dag
\\
\bMK_e &= \sin(\tfrac{u}{2}) \cos\left(\tfrac{\theta(\bNp)}{2}\right) - \epsilon\cos(\tfrac{u}{2}) \, \ba  \frac{\sin\left(\tfrac{\theta(\bN)}{2}\right)}{\sqrt{\bN}}
\end{align*}
 and with  $n\mapsto \theta(n)$  a real function such that  $\theta(0)=0$, $\forall n>0$, $\theta(n)\in]0,\pi[$  and $\lim_{n\mapsto+\infty}\theta(n)=\pi/2$.

Since we do not measure the atoms, the photon state $\rho_{k+1}$ at step $k+1$ is given by the following recurrence  from the state $\rho_k$ at step $k$:
$$
\rho_{k+1}=\bK(\rho_k)\triangleq \bM_g\rho_k\bM_g^\dag + \bM_e\rho_k\bM_e^\dag
.
$$
Consider the change of frame associated to the unitary transformation $e^{- i \hK(\bN)}$:
$
\rho = e^{- i \hK(\bN)} \rhoK e^{ i \hK(\bN)}.
$
Then we have
$
\rhoK_{k+1}=\bKK(\rhoK_{k})\triangleq\bMK_g\rhoK_k (\bMK_g)^\dag+\bMK_e\rhoK_k (\bMK_e)^\dag
.
$
It is proved in~\cite{LeghtasPhD} that, since   $|u|\leq \pi/2$,  exists a unique common eigen-state $\ket{\psiK} \in \HH_S$  of   $\bMK_g$ and $\bMK_e$.  Thus    $\rhoKb =\ket{\psiK}\bra{\psiK}$ is a  fixed point of $\bKK$. It is also proved in~\cite{LeghtasPhD} that  the  $\rhoK_{k}$'s converge to $\rhoKb$ when the function $\theta$ is strictly increasing.  Since the underlying Hilbert space $\HH_S$ is of infinite dimension, it is important to precise the type of convergence. For any initial condition $\rhoK_0$ such that  $\tr{\bN \rhoK_0} < +\infty$, then  $\lim_{k\mapsto +\infty} \tr{(\rhoK_k-\rhoKb)^2} = 0$ (Frobenius norm on Hilbert-Schmidt operators).  Since $\tr{\bN\rho} \equiv\tr{\bN\rhoK}$, we have the convergence of $\rho_k$ towards $\rho_{\infty}=e^{- i \hK(\bN)} \rhoKb e^{ i \hK(\bN)}$ as soon as the initial energy $\tr{\bN \rho_0}$ is finite:
$\lim_{k\mapsto +\infty} \tr{\left((\rho_k-\rho_\infty\right)^2} =0$. When $\theta$ is not strictly increasing, we conjecture that such convergence towards $\rho_\infty$ still holds true.

\begin{figure}
  \centering
  \includegraphics[width=0.7\textwidth]{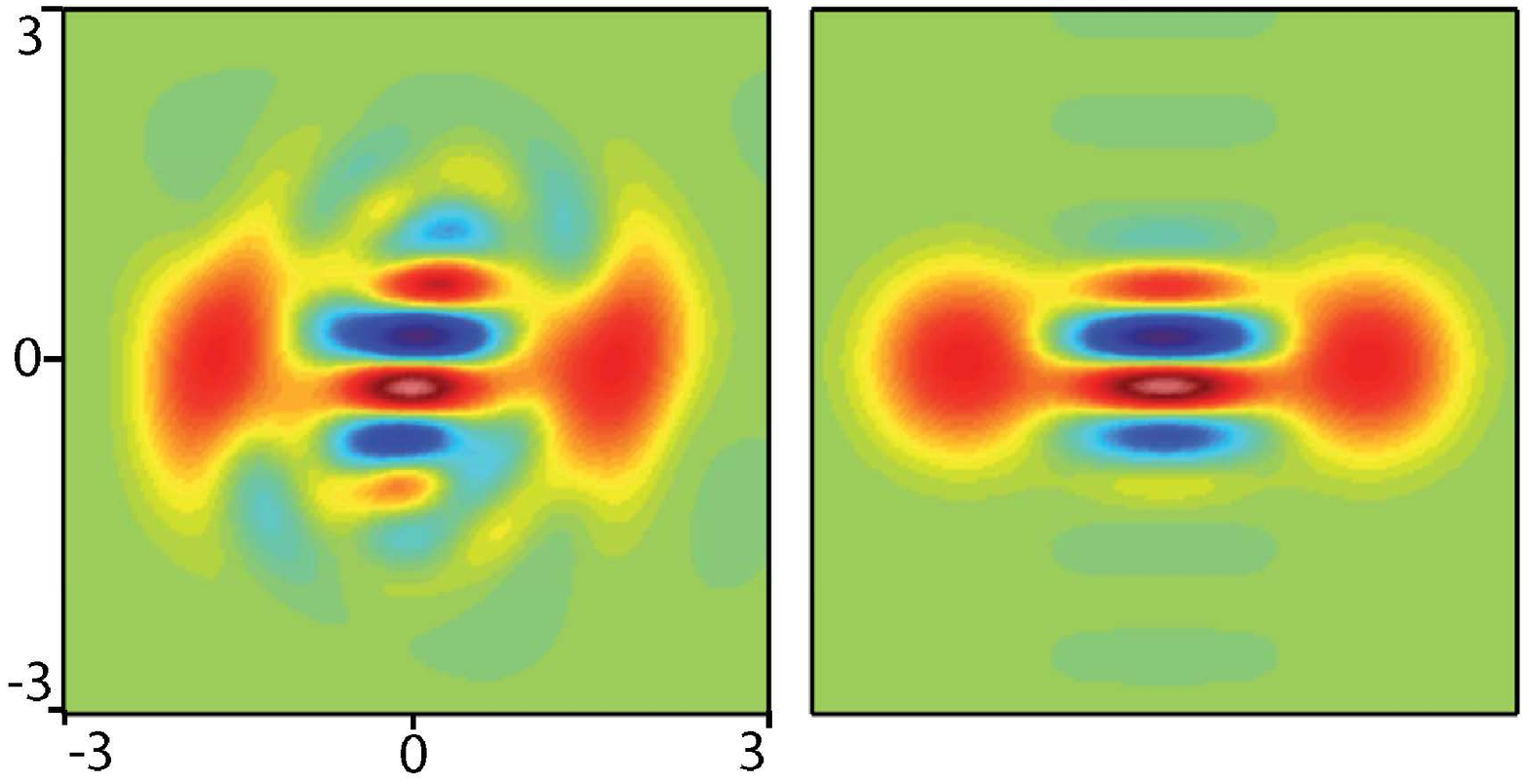}\\
  \caption{ Left: Wigner function of $\rho_\infty$ stabilized by reservoir engineering in~\cite{sarlette-et-al:PRA2012}. Right: Wigner function of a prefect Schr\"{o}dinger phase-cat,  $\tfrac{1}{2} \Big(\ket{\alpha_\infty} + i \ket{\minou\alpha_\infty}\Big) \Big( \bra{\alpha_\infty} + i \bra{\minou\alpha_\infty} \Big)$, with an  average number of photons identical to $\rho_\infty$ ($\alpha_\infty=\sqrt{\tr{\bN\rho_\infty}}$).  The color map is identical to figure~\ref{fig:Wigner_intro}.  }\label{fig:CatfullLin}
\end{figure}

For well chosen experimental  parameters~\cite{sarlette-et-al:PRA2012},   $\rhoKb$ is close to a coherent state $\ket{\alpha_\infty}\bra{\alpha_\infty}$ for some $\alpha_\infty\in\RR$ and $\hK(\bN)\approx \pi\bN^2/2$ . Since
$$
e^{-i\frac{\pi}{2} \bN^2} \ket{\alpha_\infty}=  \frac{e^{-i\pi/4}}{\sqrt{2}}\big(\ket{\alpha_\infty} + i \ket{\minou\alpha_\infty}\big),
$$ we  have under realistic conditions
$
 \lim_{k\mapsto +\infty} \rho_k\approx \tfrac{1}{2} \Big(\ket{\alpha_\infty} + i \ket{\minou\alpha_\infty}\Big) \Big( \bra{\alpha_\infty} - i \bra{\minou\alpha_\infty} \Big)
$,
 a coherent superposition of the classical states $\ket{\alpha_\infty}$ and $\ket{\minou \alpha_\infty}$ of same amplitude but of opposite phases, i.e.  a Schr\"{o}dinger phase-cat. Figure~\ref{fig:CatfullLin} displays  numerical computations of the  Wigner function of $\rho_\infty$ obtained with realistic parameters.

\section{Discrete-time  systems} \label{sec:DiscreteTime}
The theory of open quantum systems starts with the  contributions of Davies~\cite{daviesBook1976}.
The goal of this section is first  to present in an elementary way  the general structure of the  Markov models  describing such systems. Some  related stabilization problems are also addressed. Throughout this section,  $\HH$ is an Hilbert space; for each time-step $k\in\NN$, $\rho_k$ denotes  the density operator   describing the  state of  the quantum Markov process; for all $k$, $\rho_k$ is an Hilbert-Schmidt operator on $\HH$, Hermitian and of trace one; the set of continuous operators  on $\HH$ is denoted by $\LH$;  expectation values  are  denoted by the symbol $\EE{~}$.

 \subsection{Markov models}
Take a positive integer $m$ and  consider a finite set $(\bM_\mu)_{\mu \in \{1,\ldots, m\}}$ of operators on $\HH$ such that
 \begin{equation}\label{eq:PartitionUnity}
   \bI = \sum_{\mu=1}^m \bM_\mu^\dag \bM_\mu
 \end{equation}
 where $\bI$ is the identity operator.
Then each $\bM_\mu \in\LH$. Take another positive integer $m'$ and  consider a left stochastic  $m'\times m$-matrix $(\eta_{\mu'\mu}) $:  its entries are non-negative and  $\forall \mu\in\{1, \ldots, m\}$, $\sum_{\mu'=1}^{m'}\eta_{\mu'\mu} =1$. Consider the Markov process of state $\rho$ and  output $y\in\{1,\ldots, m'\}$ (measurement outcome) defined via the transition rule
\begin{equation}\label{eq:GenMarkovChain}
  \rho_{k+1} = \frac{\sum_\mu \eta_{\mu^\prime\mu} \bM_\mu \rho_k \bM_\mu^\dag}{\tr{\sum_\mu \eta_{\mu^\prime\mu} \bM_\mu \rho_k \bM_\mu^\dag}},  \quad
  y_k=\mu^\prime    \text{ with probability }  \PP_{\mu'}(\rho_k)
\end{equation}
where $\PP_{\mu'}(\rho)= \tr{\sum_\mu \eta_{\mu^\prime\mu} \bM_\mu \rho \bM_\mu^\dag}$.

\subsection{Kraus and unital maps}

The  Kraus map $\bK$  corresponds to the master equation of~\eqref{eq:GenMarkovChain}. It is given by  the expectation value of $\rho_{k+1}$ knowing $\rho_k$:
\begin{equation}\label{eq:KrausMap}
  \bK(\rho) \triangleq \sum_\mu \bM_\mu \rho \bM_\mu^\dag = \EE{\rho_{k+1}~/~\rho_k=\rho}
  .
\end{equation}
In quantum information~\cite{nielsen-chang-book}  such Kraus maps  describe quantum channels.  They admit many interesting properties. In particular,   they are contractions for many metrics (see~\cite{petz:LAA1996} for the   characterization, in finite dimension,  of  metrics for which any Kraus map is a  contraction). We just recall below two such metrics.
For any density operators $\rho$ and $\sigma$ we have
\begin{equation}\label{eq:ContractionKraus}
  D(\bK(\rho),\bK(\sigma)) \leq D(\rho,\sigma) \text{ and }  F(\bK(\rho),\bK(\sigma)) \geq F(\rho,\sigma)
\end{equation}
where the trace distance $D$ and fidelity $F$ are given by
\begin{equation}\label{eq:TraceFidelity}
  D(\rho,\sigma)\triangleq\tr{|\rho-\sigma|} \text{ and } F(\rho,\sigma) \triangleq \trr{\sqrt{\sqrt{\rho}\sigma\sqrt{\rho}}}.
\end{equation}
Fidelity  is between $0$ and $1$: $F(\rho,\sigma)=1$  if and only if, $\rho=\sigma$. Moreover $F(\rho,\sigma)=F(\sigma,\rho)$. If $\sigma=\ket\psi\bra\psi$ is a pure state ($\ket\psi$ element of $\HH$ of length one), $F(\rho,\sigma)$ coincides with the Frobenius product:
$
F(\rho, \ket\psi\bra\psi) \equiv  \tr{\rho \ket\psi\bra\psi}=\bra\psi\rho \ket\psi.
$
Kraus maps provide the evolution of  open quantum systems from an initial state $\rho_0$  without   information coming from the measurements (see~\cite[chapter 4: the environment is watching]{haroche-raimondBook06}):
$$
\rho_{k+1}=\bK(\rho_k) \text{ for } k=0,1,\ldots,
.
$$
This corresponds to the "Schr\"{o}dinger description" of the dynamics.

The "Heisenberg description" is given by the dual  map $\bK^*$. It is  characterized by $\tr{A \bK(\rho)} = \tr{\bK^*(A)\rho}$ and defined  for any operator $A$ on $\HH$  by
$$
\bK^*(A)= \sum_\mu \bM_\mu^\dag A \bM_\mu
.
$$
Technical conditions    on $A$  are required when $\HH$ is of infinite dimension, they are not given here (see, e.g., \cite{daviesBook1976}).
The map $\bK^*$ is  unital since~\eqref{eq:PartitionUnity} reads $\bK^*(\bI)=\bI$. As $\bK$, the dual map  $\bK^*$ admits a lot of interesting properties. It is noticed  in~\cite{SepulSR2010} that, based on a theorem due  of Birkhoff~\cite{Birkhoff1957},  such unital maps are contractions  on the cone of non-negative Hermitian operators equipped with the Hilbert's projective metric. In particular, when $\HH$ is of finite dimension,
 we have, for any  Hermitian  operator $A$:
 $$
 \lambda_{min}(A) \leq \lambda_{min}(\bK^*(A)) \leq \lambda_{max}(\bK^*(A)) \leq  \lambda_{max}(A)
 $$
 where $\lambda_{min}$ and $\lambda_{max}$ correspond to the smallest and largest eigenvalues.  As shown in~\cite{ReebKastoryanoWolfJMP2011},  such contraction properties  based on Hilbert's projective metric have   important implications in quantum information theory.

 To emphasize  the difference between the "Schr\"{o}dinger description" and the 'Heisenberg description" of the dynamics,  let us translate  convergence issues from the "Schr\"{o}dinger description"  to the "Heisenberg one". Assume, for clarity's sake,  that $\HH$ is of finite dimension. Suppose also  that $\bK$ admits the density operator $\bar\rho$  as unique fixed point and that, for any  initial density operator $\rho_0$,  the density operator at step $k$, $\rho_k$, defined by  $k$ iterations of $\bK$,  converges  towards $\bar\rho$  when $k$ tends to $\infty$. Then  $k\mapsto D(\rho_k,\bar\rho)$ is decreasing and converges to $0$  whereas  $k\mapsto F(\rho_k,\bar\rho)$ is increasing and converges to  $1$.

The translation of  this  convergence in the "Heisenberg description" is the following:  for any initial operator $A_0$, its $k$ iterates via $\bK^*$,  $A_k$,   converge towards $\tr{A_0\bar\rho} \bI$. Moreover when $A_0$ is Hermitian, $k\mapsto \lambda_{min}(A_k)$ and $k\mapsto \lambda_{max}(A_k)$  are  respectively increasing and decreasing  and both converge to $\tr{A_0\bar\rho}$.

\subsection{Quantum filtering}\label{ssec:filtering}

Quantum filtering has its origin in Belavkin's  work \cite{Belavkin1992} on  continuous-time  open quantum systems (see section~\ref{sec:ContinuousTime}).
The state $\rho_k$ of~\eqref{eq:GenMarkovChain} is not directly measured: open quantum systems are governed by hidden-state Markov model.
 Quantum filtering provides  an  estimate $\rhoe_k$ of $\rho_k$  based  on an initial guess  $\rhoe_{0}$ (possibly different from  $\rho_0$)  and the measurement outcomes $y_l$  between $0$ and $k-1$:
\begin{equation}\label{eq:DiscreteFilter}
  \rhoe_{l+1} = \frac{\sum_\mu \eta_{y_l\mu} \bM_\mu \rhoe_l \bM_\mu^\dag}{\tr{\sum_\mu \eta_{y_l\mu} \bM_\mu \rhoe_l \bM_\mu^\dag}},\quad  l\in\{0, \ldots, k-1\}
  .
\end{equation}
Thus   $(\rho,\rhoe)$  is  the state of  an extended  Markov process governed by  the following  rule
$$
  \rho_{k+1} = \frac{\sum_\mu \eta_{\mu^\prime\mu} \bM_\mu \rho_k \bM_\mu^\dag}{\tr{\sum_\mu \eta_{\mu^\prime\mu} \bM_\mu \rho_k \bM_\mu^\dag}}\text{ and }
     \rhoe_{k+1} = \frac{\sum_\mu \eta_{\mu^\prime\mu} \bM_\mu \rhoe_k \bM_\mu^\dag}{\tr{\sum_\mu \eta_{\mu^\prime\mu} \bM_\mu \rhoe_k \bM_\mu^\dag}}
$$
with transition probability $\PP_{\mu'}(\rho_k)= \tr{\sum_\mu \eta_{\mu^\prime\mu} \bM_\mu \rho_k \bM_\mu^\dag}$ depending  on~$\rho_k$ and independent of $\rhoe_k$.

When $\HH$ is of finite dimension, it is shown in~\cite{somaraju-et-al:acc2012} with an inequality  proved in~\cite{Rouch2011ACITo} that such discrete-time quantum filters are always stable in the following sense: the fidelity between $\rho$ and its estimate $\rhoe$ is a sub-martingale for any initial condition $\rho_0$ and $\rhoe_0$:
$
\EE{F(\rho_{k+1},\rhoe_{k+1})~|~(\rho_k,\rhoe_k)} \geq F(\rho_k,\rhoe_k)
.
$
This result  does not guaranty that $\rhoe_k$ converges to $\rho_k$ when $k$ tends to infinity. The convergence characterization of $\rhoe$ towards $\rho$  via checkable conditions on the left stochastic matrix $(\eta_{\mu'\mu})$ and on the set of operators $(\bM_\mu)$ remains  an open problem~\cite{handel:thesis,van-handel:proba2009}.

\protect{~}\\
\subsection{Stabilization via measurement-based feedback}

Assume now that the operators $\bM_\mu$  appearing in~\eqref{eq:GenMarkovChain} and satisfying~\eqref{eq:PartitionUnity},  depend  also on a control input $u$ belonging to some admissible set $\mathcal{U}$ (typically a discrete set or  a compact subset of $\RR^p$ for some positive integer $p$). Then we have the following control Markov model with input $u\in\mathcal{U}$, hidden state $\rho$ and measured output $y\in\{1,\ldots,m'\}$:
{\small\begin{equation}\label{eq:GenControlMarkovChain}
  \rho_{k+1} = \frac{\sum_\mu \eta_{\mu^\prime\mu} \bM_\mu(u_k) \rho_k \bM_\mu^\dag(u_k)}{\tr{\sum_\mu \eta_{\mu^\prime\mu} \bM_\mu(u_k) \rho_k \bM_\mu^\dag(u_k)}},  ~
  y_k=\mu^\prime    \text{ with probability }  \PP_{\mu'}(\rho_k,u_k)
\end{equation}}
where $\PP_{\mu'}(\rho,u)= \tr{\sum_\mu \eta_{\mu^\prime\mu} \bM_\mu(u) \rho \bM_\mu^\dag(u)}$.
Assume that for some nominal admissible input $\bar u\in\mathcal{U}$, this Markov process admits a steady state $\bar\rho$. This means that, for any $\mu'\in\{1,\ldots,m'\}$ we have
$
\sum_\mu \eta_{\mu^\prime\mu} \bM_\mu(\bar u) \bar\rho \bM_\mu^\dag(\bar u)= \PP_{\mu'}(\bar\rho,\bar u) \bar\rho
.
$
The  measurement-based feedback  stabilization of the steady-state $\bar\rho$  is the following problem: for any  initial condition $\rho_0$,  find  for any $k\in\NN$ a control input $u_k\in\mathcal{U}$ depending only on $\rho_0$ and on  the past $y$ values,  $(y_0,\ldots,y_{k-1})$, such  that $\rho_k$ converges almost surely towards $\bar\rho$.

Quantum-state feedback  scheme, $u=f(\rho)$, can be used here. They  can be based on Lyapunov techniques. Potential candidates  of Lyapunov functions $V(\rho)$  could be related to the metrics  for which the open-loop Kaus map with  $\bar u$ is contracting. Specific  $V$ depending on the precise structure of the system could be more adapted  as for the LKB photon box~\cite{AminiSDSMR2013A}. Such Lyapunov  feedback laws   are then given by the minimization  versus $u\in\mathcal{U}$ of
$\EE{V(\rho_{k+1 })~|~\rho_k=\rho, u_k=u}$.

Assume that  we have   a  stabilizing feedback law   $u=f(\rho)$:  $\bar u =f(\bar\rho)$ and the trajectories of \eqref{eq:GenControlMarkovChain} with $u_k=f(\rho_k)$ converge almost surely towards $\bar\rho$. Since $\rho$ is not directly accessible, one has to replace $\rho_k$ by its estimate $\rhoe_k$  to obtain $u_k$. Experimental implementations  of such  quantum feedback laws admit necessarily  an observer/controller structure governed by a  Markov process of state $(\rho,\rhoe)$ with the following transition rule:
\begin{equation}\label{eq:ObsController}
  \begin{split}
  \rho_{k+1} &= \frac{\sum_\mu \eta_{\mu^\prime\mu} \bM_\mu(f(\rhoe_k)) \rho_k \bM_\mu^\dag(f(\rhoe_k))}{\tr{\sum_\mu \eta_{\mu^\prime\mu} \bM_\mu(f(\rhoe_k)) \rho_k \bM_\mu^\dag (f(\rhoe_k))}} \\
  \rhoe_{k+1} &= \frac{\sum_\mu \eta_{\mu^\prime\mu} \bM_\mu(f(\rhoe_k)) \rhoe_k \bM_\mu^\dag(f(\rhoe_k))}{\tr{\sum_\mu \eta_{\mu^\prime\mu} \bM_\mu(f(\rhoe_k)) \rhoe_k \bM_\mu^\dag(f(\rhoe_k))}}
  \end{split}
\end{equation}
with probability  $\PP_{\mu'}(\rho_k,f(\rhoe_k))= \tr{\sum_\mu \eta_{\mu^\prime\mu} \bM_\mu(f(\rhoe_k)) \rho_k \bM_\mu^\dag(f(\rhoe_k))}$ depending on $\rho_k$ and $\rhoe_k$.   In~\cite{bouten-handel:2008} a separation principle  is proved with elementary arguments (see also~\cite{AminiSDSMR2013A}): if $\HH$ is of finite dimension, if  $\bar\rho$ is a pure state ($\bar\rho =\ket{\bar\psi}\bra{\bar\psi}$ for some $\ket{\bar\psi}$ in $\HH$)  and if $\Ker(\rhoe_0) \subset \Ker(\rho_0)$,  then  almost all realizations of~\eqref{eq:ObsController} converge to the  steady-state $(\bar\rho,\bar\rho)$.  The stabilizing feedback schemes  used in
experiments~\cite{sayrin-et-al:nature2011} and~\cite{zhouPRL2012} exploit such  observer/controller structure and rely on this separation principle where the design of the stabilizing feedback (controller)  and of the quantum-state filter (observer) are be done separately.

With such feedback scheme we loose the linear formulation of the ensemble-average master equation  with a Kraus map. In general,  there is no  simple formulation of the master equation governing  the expectation value of $\rho_k$ in closed-loop.  Nevertheless, for systems where  the  measurement step  producing the output  $y_k$ is followed by  a  control action  characterized by  $u_k$,  it is possible via a static  output  feedback, $u_k=f(y_k)$ where $f$ is now  some function from $\{1,\ldots,m'\}$ to $\mathcal{U}$, to  preserve  in closed-loop such  Kraus-map formulations. These specific  feedback  schemes, called Markovian feedbacks,   are due to Wiseman and have important applications. They are well explained and illustrated in the  recent book~\cite{wiseman-milburnBook}.

\subsection{Stabilization of pure states  by reservoir engineering}
With  $T$ as  sampling period, a possible formalization of  this passive  stabilization method is as  follows.
The goal is to stabilize a  pure  state  $\bar\rho_S=\ket{\bar\psi_S}\bra{\bar\psi_S}$  for a system $S$ with Hilbert space $\HH_S$  and Hamiltonian operator  $H_S$ ($\ket{\bar\psi_S}\in\HH_S$ is of length one).  To achieve this goal  consider   a "realistic"  quantum controller   of Hilbert space $\HH_C$ with initial state $\ket{\theta_C}$ and  with Hamiltonian $\bH_C$. One has  to design  an adapted  interaction between $S$ and $C$  with a well chosen interaction Hamiltonian  $\bH_{int}$, an Hermitian operator on $\HH_{S,C}=\HH_S\otimes \HH_C$.
This controller  $C$ and its interaction with $S$ during the sampling interval  of length $T$ have to fulfill the  conditions explained  below  in order to stabilize $\bar\rho_S$.

 Denote by  $\bU_{S,C}=\bU(T)$ the  propagator between $0$ and time $T$ for  the composite system $(S,C)$: $\bU(t)$ is the unitary operator on $\HH_{S,C}$ defined  by
 $$
 \dotex\bU = - \tfrac{i}{\hbar} \Big(\bH_S\otimes \bI_C+\bH_{int}+\bI_S\otimes \bH_C\Big) \bU, \quad \bU(0)=\bI_{S,C}
 $$
 where $\bI_S$, $\bI_C$ and $\bI_{S,C}$ are the identity operators on $\HH_S$, $\HH_C$,  and $\HH_{S,C}$, respectively.
To the propagator $U_{S,C}$ and the initial controller wave function $\ket{\theta_C}\in\HH_C$ is attached a   Kraus map $\bK$ on $\HH_S$,
$$
   \bK(\rho_S)= \sum_\mu \bM_\mu \rho_S \bM_\mu^\dag
$$
where the   operators $\bM_\mu$ on $\HH_S$   are defined  by the decomposition,
$$
\forall \ket{\psi_S}\in\HH_S, \quad \bU_{S,C}\big( \ket{\psi_S}\otimes\ket{\theta_C}\big)= \sum_\mu \big(\bM_{\mu} \ket{\psi_S}\big)\otimes \ket{\lambda_\mu}
,
$$
with  $(\ket{\lambda_{\mu}})$ any  ortho-normal basis of  $\HH_C$. Despite the fact that the operators $(\bM_\mu)$ depend  on the choice of  this basis, the Kraus map   $\bK$ is independent of  this choice: it depends only on $U_{S,C}$ and $\ket{\theta_C}$.

The first stabilization condition is the following:   the Kraus  operators $\bM_\mu$  have to admit  $\ket{\bar\psi_S}$ as a common  eigen-vector since $\bar\rho_S$ has to be  a fixed point of $\bK$ ($\bK(\bar\rho_S)=\bar\rho_S$).

The second stabilization condition  is the following:  for   any initial  density operator $\rho_{S,0}$,  the  iterates $\rho_{S,k}$  of $\bK$ converge to $\bar\rho_S$, i.e.,
$$
\lim_{k\mapsto +\infty} \rho_{S,k}=\bar\rho_S \text{ where } \rho_{S,k}=\bK(\rho_{S,k-1})
.
$$

 When these two conditions are satisfied, the repetition of   the same  interaction for  each sampling interval $[kT,(k+1)T]$ ($k\in\NN)$  with a  controller-state $\ket{\theta_C}$ at  $kT$ ensures that the density operator of  $S$ at  $kT$, $\rho_{S,k}$,  converges to $\bar\rho_S$ since $\rho_{S,k}=\bK(\rho_{S,k-1})$. Here,  the so-called  reservoir is made of the infinite set  of identical  controller systems $C$ indexed by $k\in\NN$,  with initial state $\ket{\theta_C}$ and interacting sequentially  with $S$ during $[kT,(k+1)T]$.

\section{Continuous-time  systems}  \label{sec:ContinuousTime}

\subsection{Stochastic master equations}
These models have their origins in the work of Davies~\cite{daviesBook1976}, are related to quantum trajectories~\cite{carmichael-book,dalibard-et-al:PRL92} and are connected to  Belavkin  quantum filters~\cite{Belavkin1992}. A modern and  mathematical  exposure of the diffusive models   is given in~\cite{BarchielliGregorattiBook}.
These models are interpreted  here  as  continuous-time versions  of~\eqref{eq:GenMarkovChain}. They are based on    stochastic differential equations, also called Stochastic Master Equations (SME). They provide  the evolution of the  density operator $\rho_t$ with respect to the time $t$. They are driven by a finite number of independent  Wiener processes indexed by $\nu$,   $(W_{\nu,t})$, each of them being  associated to a continuous classical and real  signal, $y_{\nu,t}$, produced  by detector $\nu$.  These SMEs admit the following form:
\begin{multline} \label{eq:DiffusiveSME}
d\rho_t=\left(-\tfrac{i}{\hbar}[\bH,\rho_t]+\sum_\nu \bL_\nu \rho_t \bL_\nu^\dag - \tfrac{1}{2} (\bL_\nu^\dag \bL_\nu\rho_t+\rho_t \bL_\nu^\dag \bL_\nu)\right) dt
\\
+ \sum_\nu \sqrt{\eta_\nu}\bigg(\bL_\nu\rho_t+\rho_t \bL_\nu^\dag-\tr{(\bL_\nu+\bL_\nu^\dag)\rho_t}\rho_t\bigg) dW_{\nu,t}
\end{multline}
where $\bH$ is the Hamiltonian operator on the underlying Hilbert space $\HH$ and $\bL_\nu$ are arbitrary operators  (not necessarily Hermitian) on $\HH$.  Each measured signal $y_{\nu,t}$ is related to $\rho_t$ and $W_{\nu,t}$  by the following  output relationship:
$$
dy_{\nu,t}= dW_{\nu,t}+\sqrt{\eta_\nu}\tr{(\bL_\nu+\bL_\nu^\dag)\,\rho_t}\,dt
$$
where  $\eta_\nu \in[0,1]$ is  the  efficiency of detector  $\nu$.
The ensemble average of $\rho_t$ obeys  thus to a linear differential  equation, also called master  or   Lindblad-Kossakowski  differential equation~\cite{Kossakowski-1972,Lindblad-1976}:
\begin{equation}\label{eq:Lindblad}
\dotex\rho=-\tfrac{i}{\hbar}[\bH,\rho]+\sum_\nu \bL_\nu \rho_t \bL_\nu^\dag - \tfrac{1}{2} (\bL_\nu^\dag \bL_\nu\rho_t+\rho_t \bL_\nu^\dag \bL_\nu)
.
\end{equation}
It is the continuous-time analogue of the Kraus map~$\bK$ associated to the Markov process~\eqref{eq:MarkovChainError}.

In fact~\eqref{eq:GenMarkovChain} and~\eqref{eq:DiffusiveSME} have the same  structure. This becomes obvious if one remarks that, with standard  It$\bar{\text{o}}$ rules, \eqref{eq:DiffusiveSME} admits the following formulation
$$
\rho_{t+dt}= \frac{
 \bM_{{\!{dy_{t}}}}\rho_t \bM_{{\!{dy_{t}}}}^\dag + \sum_\nu (1-\eta_\nu ) \bL_{\nu}\rho_t \bL_{\nu}^\dag dt
                          }{
                          \tr{
 \bM_{{\!{dy_{t}}}}\rho_t \bM_{{\!{dy_{t}}}}^\dag + \sum_\nu (1-\eta_\nu ) \bL_{\nu}\rho_t \bL_{\nu}^\dag dt
                          }}
$$
with $ \bM_{{\!{dy_{t}}}} = \bI + \left(-\tfrac{i}{\hbar}\bH -\tfrac{1}{2} \sum_\nu \bL_{\nu}^\dag \bL_{\nu}\right) dt   + \sum_\nu \sqrt{\eta_\nu} {{dy_{\nu t}}} \bL_{\nu}$.  Moreover  the  probability  associated to the measurement outcome  $dy=(dy_{\nu})$,  is given by the following density
\begin{multline*}
  \PP \left( dy  \in\prod_{\nu} [\xi_\nu,\xi_\nu+d\xi_\nu]~\Big/~\rho_t\right)
  \\
  =\tr{\bM_{{\!{\xi}}} \rho_t \bM_{{\!{\xi}}}^\dag ++ \sum_\nu (1-\eta_\nu ) \bL_{\nu}\rho_t \bL_{\nu}^\dag dt}~\prod_\nu  e^{- \xi_\nu^2/2dt}\tfrac{d\xi_\nu}{\sqrt{2\pi dt}}
\end{multline*}
where $\xi$ stands for the vector $(\xi_\nu)$.
 With such a  formulation, it becomes clear that~\eqref{eq:DiffusiveSME} preserves the  trace and the non-negativeness of $\rho$.  This formulation provides also directly  a time discretization numerical  scheme  preserving  non-negativeness of $\rho$ (see appendix~\ref{ap:num}).

Mixed diffusive/jump stochastic master equations can be considered.  Additional Poisson counting processes $(N_\mu(t))$  are  added  in parallel to  the  Wiener processes $(W_{\nu,t})$~\cite{AminiBelavkin2014}:
{\small \begin{multline} \label{eq:DiffusiveJumpSME}
d\rho_t =\left(-\tfrac{i}{\hbar}[\bH,\rho_t] +  \sum_\nu \bL_{\nu} \rho_t \bL_{\nu}^\dag - \tfrac{1}{2}(\bL_{\nu}^\dag \bL_{\nu}\rho_t+\rho_t \bL_{\nu}^\dag \bL_{\nu}) \right)\,dt
\\
+ \sum_\nu  \sqrt{\eta_\nu} \bigg(\bL_{\nu}\rho_t+\rho_t \bL_{\nu}^\dag-\tr{(\bL_{\nu}+\bL_{\nu}^\dag)\rho_t}\rho_t\bigg){{dW_{\nu,t}}}
\\
+ \left(\sum_\mu \bV_{\mu}\rho_t \bV_{\mu}^\dag - \tfrac{1}{2}(\bV_{\mu}^\dag \bV_{\mu}\rho_t +\rho_t \bV_{\mu}^\dag \bV_{\mu}) \right)\,dt
\\
 +  \sum_\mu \left(\tfrac{\thetab_\mu\rho_t+\sum_{\mu'} \etab_{\mu,\mu'} \bV_{\mu'}\rho_t \bV_{\mu'}^\dag}{\thetab_\mu +
   \sum_{\mu'}\etab_{\mu,\mu'}\tr{\bV_{\mu'}\rho_t \bV_{\mu'}^\dag}} -\rho_t\right)
      \left({{dN_\mu(t) }}-\Big( \thetab_\mu + \sum_{\mu'}\etab_{\mu,\mu'} \tr{\bV_{\mu'}\rho_t \bV_{\mu'}^\dag}\Big)\,dt\right)
\end{multline}}
where the $\bV_\mu$'s are  operators on $\HH$,  where the additional parameters   $\thetab_\mu, \etab_{\mu,\mu'}\geq 0$ with $\etab_{\mu'}=\sum_{\mu} \etab_{\mu,\mu'} \leq 1$,  describe  counting  imperfections.  For each~$\mu$,  $\Big( \thetab_\mu + \sum_{\mu'}\etab_{\mu,\mu'} \tr{\bV_{\mu'}\rho_t \bV_{\mu'}^\dag}\Big)\,dt$ is  the probability  to increment by one  $N_{\mu}$  between $t$ and $t+dt$.

For any vector $\xi=(\xi_\nu)$, take the following definition for $\bM_\xi$
$$
\bM_{{{\xi}}} = \bI - \left(\tfrac{i}{\hbar}\bH +\tfrac{1}{2} \sum_\nu \bL_{\nu}^\dag \bL_{\nu} + \tfrac{1}{2} \sum_\mu \bV_{\mu}^\dag \bV_{\mu}\right) dt + \sum_\nu \sqrt{\eta_\nu} {{\xi_{\nu}}} \bL_{\nu}
$$
and consider the following partial Kraus map  depending on $\xi$: 
$$
\bK_\xi(\rho) = \bM_{{{\xi}}}\rho \bM_{{{\xi}}}^\dag + \sum_\nu (1-\eta_\nu ) \bL_{\nu}\rho \bL_{\nu}^\dag dt+  \sum_\mu (1-\etab_\mu ) \bV_{\mu} \rho \bV_{\mu}^\dag dt
.
$$
The stochastic model~\eqref{eq:DiffusiveJumpSME}   is   similar to the discrete-time Markov process~\eqref{eq:GenMarkovChain} where the  discrete-time  outcomes $y_k$  is replaced by the continuous-time outcomes   $(dy_t,dN(t))$. More precisely,  the transition from $\rho_t$ to $\rho_{t+dt}$ is given by the following  transition rules:
\begin{enumerate}

\item  The transition  corresponding  to no-jump   outcomes $(dy_t,dN(t)=0)$ reads
$$
\rho_{t+dt}= \frac{\bK_{dy_t}(\rho_t)}{\tr{\bK_{dy_t}(\rho_t)}}
$$
 and is associated to the following probability law:
{\small
\begin{multline*}
  \PP \left( dy  \in\prod_{\nu} [\xi_\nu,\xi_\nu+d\xi_\nu]\text{ and } dN(t)=0~\bigg/~\rho_t\right)
  \\
  =\left(1-\left(\sum_\mu \thetab_\mu\right)dt\right)\tr{\bK_\xi(\rho_t)} \left(\prod_\nu  e^{- \xi_\nu^2/2dt}\tfrac{d\xi_\nu}{\sqrt{2\pi dt}}\right)
\end{multline*}
}
Since
$$
\int_\xi \tr{\bK_\xi(\rho_t)}~\prod_\nu  e^{- \xi_\nu^2/2dt}\tfrac{d\xi_\nu}{\sqrt{2\pi dt}}
= 1-  \left(\etab_{\mu}\tr{\bV_{\mu}\rho_{t} \bV_{\mu}^\dag}\right) dt + O(dt^2)
$$
 we recover the usual  no-jump probability, $1- \left(\sum_\mu \thetab_\mu+ \etab_{\mu}\tr{\bV_{\mu}\rho_{t} \bV_{\mu}^\dag} \right) dt$,  up to $O(dt^2)$ terms.

\item   The transition  corresponding  to   outcomes with  a single  jump of label $\mu$,   $(dy_t,dN(t)=(\delta_{\mu,\mu'})_{\mu'})$,  reads
$$
\rho_{t+dt}= \frac{
     \bK_{dy_t}\left( \thetab_\mu\rho_{t}+\sum_{\mu'}\etab_{\mu,\mu'} \bV_{\mu'}\rho_{t} \bV_{\mu'}^\dag \right)
                         }{
   \tr{ \bK_{dy_t}\left( \thetab_\mu\rho_{t}+\sum_{\mu'}\etab_{\mu,\mu'} \bV_{\mu'}\rho_{t} \bV_{\mu'}^\dag \right)}
                          }
$$
 and is associated to the following probability law:
{\small
\begin{multline*}
   \PP \left( dy  \in\prod_{\nu} [\xi_\nu,\xi_\nu+d\xi_\nu]\text{ and } dN(t)=(\delta_{\mu,\mu'})_{\mu'}~\bigg/~\rho_t\right)
  \\
  =dt~ \tr{ \bK_{\xi}\left( \thetab_\mu\rho_{t}+\sum_{\mu'}\etab_{\mu,\mu'} \bV_{\mu'}\rho_{t} \bV_{\mu'}^\dag \right)}\left(\prod_\nu  e^{- \xi_\nu^2/2dt}\tfrac{d\xi_\nu}{\sqrt{2\pi dt}}\right)
\end{multline*}
}
By integration versus $\xi$, we recover, up to $O(dt^2)$ terms, the probability of jump $\mu$: $\left(\thetab_\mu+\sum_{\mu'} \etab_{\mu,\mu'}\tr{\bV_{\mu'}\rho_{t+dt} \bV_{\mu'}^\dag} \right) dt$. 

\item The probability to have at least  two  jumps, i.e.   $dN_{\mu}(t)=dN_{\mu'}(t)=1$ for some  $\mu\neq \mu'$,   is an $O(dt^2)$ and thus  negligible.   

\end{enumerate}
Standard  computations show that such time  discretization schemes   converge in law to the continuous-time process~\eqref{eq:DiffusiveJumpSME} when $dt$ tends to $0$.  They preserve the fact that  $\rho\geq 0$ and can be used for Monte-Carlo simulations and quantum filtering.

\subsection{Quantum filtering}
For clarity's sake,  take in~\eqref{eq:DiffusiveSME}  a single  measurement $y_t$ associated to  operator $\bL$,  detection efficiency $\eta \in[0,1]$ and scalar Wiener process $W_t$: ${{dy_t}}=  \sqrt{\eta} \tr{(\bL+\bL^\dag)\,\rho_t}\,dt + {{dW_t}}$. The continuous-time counterpart  of~\eqref{eq:DiscreteFilter}  provides
the estimate $\rhoe_{t}$ by the   Belavkin quantum filtering process
\begin{multline*}
d\rhoe_{t} =-\tfrac{i}{\hbar}[\bH,\rhoe_{t}] \,dt+\left(\bL\rhoe_{t} \bL^\dag - \tfrac{1}{2}(\bL^\dag \bL\rhoe_{t} +\rhoe_{t} \bL^\dag \bL ) \right)\,dt
\\
+\sqrt{\eta}
 \left( \bL\rhoe_{t} + \rhoe_{t}  \bL^\dag - \tr{(\bL+\bL^\dag)\rhoe_{t}} \rho_t^e\right)
\left({{dy_t}}-   \sqrt{\eta}\tr{(\bL+\bL^\dag)\rhoe_{t}}\,dt\right).
\end{multline*}
initialized to any  density matrix $\rhoe_{0}$. Thus $(\rho,\rhoe)$ obeys to the following  set of nonlinear stochastic  differential equations
  \begin{align*}
&d\rho_t =-\tfrac{i}{\hbar}[\bH,\rho_t] \,dt+\left(\bL \rho_t \bL^\dag - \tfrac{1}{2}(\bL^\dag \bL\rho_t +\rho_t \bL^\dag \bL )\right) \,dt
\\&\quad + \sqrt{\eta}\left( \bL\rho_t + \rho_t \bL^\dag - \tr{(\bL+\bL^\dag)\rho_t} \rho_t \right) \, {{dW_t}}
\\[1.em]
&d\rhoe_{t} =-\tfrac{i}{\hbar}[\bH,\rhoe_{t}] \,dt+\left(\bL\rhoe_{t}\bL^\dag - \tfrac{1}{2}(\bL^\dag \bL\rhoe_{t} +\rhoe_{t} \bL^\dag \bL ) \right)\,dt
\\
&+\sqrt{\eta}
 \left( \bL\rhoe_{t} + \rhoe_{t} \bL^\dag - \tr{(\bL+\bL^\dag)\rhoe_{t}} \rhoe_{t}\right){{dW_t}}
 \\
 &+{\eta
 \left( \bL\rhoe_{t} + \rhoe_{t} \bL^\dag - \tr{(\bL+\bL^\dag)\rhoe_{t}} \rhoe_{t}\right) \tr{(\bL+\bL^\dag){(\rho_t-\rhoe_{t})}}\,dt }
 .
\end{align*}
 It is proved in~\cite{AminiBelavkin2014} that  such filtering process is always stable in the sense that, as for the discrete-time case, the fidelity between $\rho_t$ and $\rhoe_t$ is a sub-martingale. In~\cite{van-handel:proba2009} a first convergence analysis of these filters is proposed. Nevertheless  the convergence characterization in  terms of the operators $\bH$, $\bL$ and  the parameter $\eta$ remains an open problem as far as we know.

 Formulations of quantum filters  for   stochastic master equations   driven by  an arbitrary number of Wiener and Poisson  processes can be found in~\cite{AminiBelavkin2014}.

\subsection{Stabilization via measurement-based feedback}

Assume that the Hamiltonian $H=H_0+ u H_1 $ appearing in~\eqref{eq:Lindblad} depends on some scalar control input   $u$, $H_0$ and $H_1$ being Hermitian operators on $\HH$. Assume also that $\bar\rho=\ket{\bar\psi}\bra{\bar\psi}$ is a steady-state of~\eqref{eq:Lindblad} for $u=0$.  Necessarily  $\ket{\bar\psi}$ is an eigen-vector of each $\bL_\nu$, $\bL_\nu\ket{\bar\psi}=\lambda_\nu \ket{\bar\psi}$ for some $\lambda_\nu \in\CC$.  This implies that $\bar\rho$ is also a steady-state of~\eqref{eq:DiffusiveSME} with $u=0$, since $\bL_\nu\bar\rho + \bar\rho \bL_\nu^\dag = \tr{(\bL_\nu+\bL_\nu^\dag)\bar\rho} \bar\rho$. The stabilization of $\bar\rho$ consists then in finding a feedback law $u=f(\rho)$ with $f(\bar\rho)=0$ such that almost all trajectories $\rho_t$ of the closed-loop  system~\eqref{eq:DiffusiveSME} with $H=H(t)=H_0+f(\rho_t) H_1$  converge to $\bar\rho$ when $t$ tends to $+\infty$. Such feedback law could be obtained  by Lyapunov techniques as in~\cite{mirrahimi-handel:siam07}. As in the discrete-case, $\rho_t$ is replaced, in the feedback law,  by its estimate $\rhoe_t$ obtained via quantum filtering. Convergence  is then guarantied as soon as $\Ker{\rhoe_0} \subset \Ker{\rho_0}$ \cite{bouten-handel:2008}.
Other  feedback schemes not relying   directly on the  quantum state $\rho_t$  but still  based on past values of the measurement signals $y_{\nu}$ can be considered (see~\cite{wiseman-milburnBook} for  Markovian feedbacks; see~\cite{VijayMSWMNKS2012N,CampaFRDMMDMH2013PRX} for recent experimental implementations).

\subsection{Stabilization via coherent feedback}

This passive  stabilization method has its origin, for classical system,  in the classical  Watt regulator   where a mechanical system,  the  steam machine,  was controlled by another  mechanical system,  a conical pendulum. As initially shown in~\cite{maxwell-1868}, the study of  such closed-loop systems  highlights stability and convergence as the main mathematical issues.  For quantum systems, these issues remain similar  and are   related to reservoir engineering~\cite{PoyatCZ1996PRL,Lloyd2000PRA}.

As in the discrete-time case, the goal remains to stabilize a  pure  state  $\bar\rho_S=\ket{\bar\psi_S}\bra{\bar\psi_S}$  for system  $S$  (Hilbert space $\HH_S$  and Hamiltonian   $\bH_S$) by coupling  to the controller system $C$ (Hilbert space $\HH_C$, Hamiltonian $\bH_C$)  via the interaction $\bH_{int}$, an Hermitian operator on $\HH_S\otimes \HH_C$.  The controller $C$ is subject to decoherence described  by the  set $\Big(\bL_{C,\nu}\Big)$ of operators on $\HH_C$ indexed by $\nu$. The closed-loop system is a composite system with Hilbert space $\HH=\HH_S\otimes\HH_C$. Its density operator $\rho$ obeys to~\eqref{eq:Lindblad} with
$\bH=\bH_S\otimes \bI_C+\bI_S\otimes \bH_C + \bH_{int}$ and $\bL_\nu = \bI_S\otimes\bL_{C,\nu}$
( $\bI_S$ and $\bI_C$  identity operators on $\HH_S$ and $\HH_C$,  respectively).  Stabilization is achieved when $\rho(t)$ converges, whatever its initial condition $\rho(0)$ is, to a separable state of the form $\bar\rho_S \otimes \bar\rho_C$ where $\bar\rho_C$  could possibly depend on $t$  and/or on~$\rho(0)$.  In several interesting cases, such as  cooling~\cite{HamerM2012PRL},  coherent feedback is  shown to  outperform measurement-based feedback.

The asymptotic analysis (stability and  convergence rates)  for  such composite  closed-loop systems is   far from being obvious, even if such analysis  is based  on known properties for  each subsystem and for  the coupling Hamiltonian $\bH_{int}$.

When $\HH$ is of infinite dimension,  convergence   analysis becomes more difficult. To have an idea of the  mathematical issues,  we will consider two examples of physical interest.  The first one is derived form~\cite{sarlette-et-al:PRA2012}:
\begin{equation}\label{eq:CatKerr}
	\tfrac{d}{dt}\rho =u [ \ba^\dag-\ba, \rho]  +\kappa \big(\ba \rho \ba^\dag - (\bN \rho+\rho \bN)/2  \big)
  + \kappa_c \big( e^{i\pi\bN}\ba \rho \ba^\dag e^{-i\pi\bN}  - (\bN \rho +\rho \bN)/2   \big)
\end{equation}
where $u$, $\kappa$ and $\kappa_c$ are strictly positive parameters. It is shown in~\cite{sarlette-et-al:PRA2012}, that~\eqref{eq:CatKerr} admits a unique   steady state $\rho_{\infty}$ given  by  its Glauber-Shudarshan $P$ distribution:
$$
\rho_{\infty} = \int_{-\alpha^c_\infty}^{\alpha^c_\infty} \mu(x) \ket{x}\bra{x}~ dx
$$
where $\ket{x}$ is the coherent state of real amplitude $x$ and where the  non-negative  weight function $\mu$   reads
$$
\mu(x) = \mu_0\, \frac{\left(((\alpha^c_\infty) ^2 - x^2)^{(\alpha^c_\infty)^2}~e^{ x^2}\right)^{r_c}}{\alpha^c_\infty  - x}\ ,
$$
with $r_c=2\kappa_c/(\kappa+\kappa_c)$ and $\alpha^c_\infty=2u/(\kappa+\kappa_c)$.
The normalization factor $\mu_0>0$ ensures that $ \int_{-\alpha^c_\infty}^{\alpha^c_\infty} \mu(x) dx\, =1$, i.e., $\tr{\rho_\infty}=1$.
We conjecture that any solution $\rho(t)$  of~\eqref{eq:CatKerr} starting from any initial condition $\rho(0)=\rho_0$ of finite energy ($\tr{\rho_0 \bN} < \infty$), converges in Frobenius norm towards   $\rho_{\infty}$.  When $\rho$ follows~\eqref{eq:CatKerr}, its  Wigner function $W^\rho$ (see appendix~\ref{ap:oscillator}) obeys to the following  Fokker-Planck equation  with non-local terms ($\Delta= \dvv{~}{x}+\dvv{~}{p}$):
\begin{multline*}
   \left. \dv{W^{\rho}}{t}\right|_{(t,x,p)} =
    \tfrac{\kappa+\kappa_c}{2} \left(\dv{~}{x} \Big( (x-\alpha_\infty) W^{\rho}\Big) +  \dv{~}{p} \Big( p W^{\rho}\Big)
    + \tfrac{1}{4}  \Delta W^{\rho} \right)_{(t,x,p)}
    \\
+ \kappa_c \left((x^2+p^2+\tfrac{1}{2})\left( \left.W^{\rho}\right|_{{{(t,-x,-p)}}} -  \left.W^{\rho}\right|_{(t,x,p)} \right) + \tfrac{1}{16}\left( \left.\Delta W^{\rho}\right|_{{{(t,-x,-p)}}} - \left.\Delta W^{\rho} \right|_{(t,x,p)} \right)
    \right)
    \\
 - \kappa_c \left(  \tfrac{x}{2}\left( \left.\dv{W^{\rho}}{x}\right|_{{{(t,-x,-p)}}}+ \left.\dv{W^{\rho}}{x}\right|_{(t,x,p)} \right)
        +  \tfrac{p}{2}\left( \left.\dv{W^{\rho}}{p}\right|_{{{(t,-x,-p)}}}+ \left.\dv{W^{\rho}}{p}\right|_{(t,x,p)} \right)  \right)
        .
\end{multline*}
This partial differential equation is derived from the correspondence relationships~\eqref{eq:correspondence} and   $W^{e^{i\pi\bN}\rho e^{-i\pi\bN}  } (x,p)\equiv W^{\rho}(-x,-p)$.  We conjecture that $W^\rho(t,x,p)$ converges, when  $t \mapsto +\infty$,
 towards $$W^{\rho_\infty}(x,p) =  \int_{-\alpha^c_\infty}^{\alpha^c_\infty} \tfrac{2 \mu(\alpha)}{\pi} e^{-2 (x-\alpha)^2-2p^2}~ d\alpha$$ for any initial condition $W_0=W^{\rho_0}$ with finite energy, i.e., such that  (see, e.g.,\cite{haroche-raimondBook06}[equation (A.42)]),
$$
 \iint_{\RR^2}  (x^2+p^2)W_0(x,p) ~dx dp = \tfrac{1}{2}+\tr{\bN\rho_0} < +\infty
 .
$$

The second example is derived  from~\cite{MirrahimiCatComp2014} and  could have important applications for quantum computations. It is  governed by the following  master equation:
\begin{equation}\label{eq:LindbladCatQubit}
 \dotex \rho = u [ (\ba^\dag)^r - \ba^r,\rho] + \kappa \left(( \ba^r \rho (\ba^\dag)^r - \tfrac{1}{2}(\ba^\dag)^r\ba^r \rho -  \tfrac{1}{2}\rho (\ba^\dag)^r\ba^r \right)
\end{equation}
where $u> 0$ and $\kappa>0$ are constant parameters  and $r$ is an integer greater than~$1$. Set $\bar\alpha= \sqrt[r]{2u/\kappa}$ and for $s\in\{0,1,\ldots,r-1\}$,
$\bar\alpha_s=e^{2is\pi/r} \bar\alpha$. Denote by $\ket{\bar\alpha_s}$ the coherent state of  complex  amplitude $\bar\alpha_s$.
Computations exploiting  properties of coherent states recalled  in appendix~\ref{ap:oscillator}  show  that, for any  $s$,  $\ket{\bar\alpha_s}\bra{\bar\alpha_s}$ is a steady state of~\eqref{eq:LindbladCatQubit}. Moreover the set of steady states corresponds to   the density operators $\bar\rho$ with support inside the vector space spanned by  the $\ket{\bar\alpha_s}$ for $s\in\{0,1,\ldots,r-1\}$. We  conjecture that, for  initial conditions $\rho(0)$  with finite energy ($\tr{\rho \bN} <\infty$),  the solutions of~\eqref{eq:LindbladCatQubit} are well defined and  converge in Frobenius norm to   such steady states  $\bar\rho$  possibly depending on $\rho(0)$. Having sharp  estimations of the convergence rates is also an open question.  We cannot apply here the existing  general convergence results towards  "full  rank  steady-states"  (see, e.g., \cite{AttalBook3_2006}[chapter 4]): here  the rank of such steady states $\bar\rho$ is at most $r$.
Another formulation of such dynamics  can be  given via   the Wigner  function $W^\rho$  of   $\rho$ (see appendix~\ref{ap:oscillator}). With  the correspondence~\eqref{eq:correspondence}, \eqref{eq:LindbladCatQubit} yields  a partial differential equation   describing the time evolution of $W^\rho$: this equation is of order one in time but of order $2r$ versus the phase plane variables $(x,p)$. It corresponds to an unusual  Fokker-Planck equation of high order.

\section{Concluding remarks}

The above exposure deals  with specific and limited aspects  of modelling and control of open quantum systems.   It does not  consider many other interesting developments such as
\begin{itemize}
  \item controllability and motion planing   in finite dimension~\cite{alessandro:book,grigoriu:hal-00696546}  and in  infinite dimension (see, e.g., \cite{beauchard-coron:06,beauchard-et-al:JMP2010,BeaucPR2013A,chambrion-et-al:IHP09,ervedoza-puel:ihp2009});
  \item quantum Langevin equations and input/output approach~\cite{gardiner-zollerBook}, quantum signal amplification~\cite{ClerkDGMS2010RMP}  and linear quantum systems~\cite{JamesNP2008ACITo};

  \item  $(S,L,H)$ formalism for quantum networks~\cite{GoughJ2009ACITo};

  \item  master equations and quantum Fokker Planck equations~\cite{CarmichaelBook99,CarmichaelBook07};
  \item optimal control methods~\cite{NielsKGK2010,baudoin-salomon:scl2010, baudoinJDE2005,BonnaCGLSZ2012ACITo,GaronGS2013PRA}.
\end{itemize}
More topics can also be found in the   review articles~\cite{RNC:RNC1016,James2011,AltafT2012ACITo}.

\appendix

\section{Quantum harmonic oscillator}\label{ap:oscillator}

 We just recall here some useful formulae (see, e.g.,~\cite{barnett-radmoreBook}). The Hamiltonian formulation of the classical harmonic oscillator of pulsation  $\omega >0$, $\ddotex x = - \omega^2 x $, is as follows:
$$
\dotex x =\omega p= \dv{H}{p}, \quad \dotex p = -\omega x = -
\dv{H}{x} \quad
$$
with the classical Hamiltonian $H(x,p) =\frac{\omega}{2} (p^2 + x^2)$.  The correspondence principle yields the following quantization:
$H$ becomes an operator $\bH$ on the function of $x\in\RR$   with complex values. The classical state $(x(t),p(t))$ is replaced by  the quantum state $\ket{\psi}_t$ associated to the  function   $\psi(x,t)\in\CC$. At each  $t$, $\RR\ni x \mapsto \psi(x,t)$ is measurable and  $\int_{\RR} |\psi(x,t)|^2 dx =1$: for each $t$,  $\ket{\psi}_t \in L^2(\RR,\CC)$.

The Hamiltonian  $\bH$   is derived from  the classical one   $H$ by replacing  $x$ by   the Hermitian operator   $\bX\equiv \tfrac{x}{\sqrt{2}}$ and  $p$ by the Hermitian operator    $\bP\equiv-\tfrac{ i }{\sqrt{2}}  \dv{}{x}$:
 $$
 \frac{\bH}{\hbar}=\omega (\bP^2 + \bX^2) \equiv -\frac{\omega}{2}\dvv{}{x} +
 \frac{\omega}{2} x^2
 .
$$
The Hamilton  ordinary differential equations are replaced by the  Schr\"{o}dinger equation,
$\dotex \ket{{\psi}} = - \imath\tfrac{\bH}{\hbar} \ket{{\psi}}$,
a partial differential equation defining $\psi(x,t)$  from its initial condition  $(\psi(x,0))_{x\in\RR}$:
$
 \imath   \dv{\psi}{t}(x,t)  = - \frac{\omega}{2} \dvv{\psi}{x}(x,t) +
\frac{\omega}{2} x^2 \psi(x,t),\quad x\in\RR
.
$
The average position reads
$
 \bket{\bX}_t = \bket{\psi|\bX|\psi} = \tfrac{1}{\sqrt{2}}\int_{-\infty}^{+\infty} x|\psi|^2 dx.
$
The average impulsion reads
$
 \bket{\bP}_t = \bket{\psi|\bP|\psi}= -\tfrac{ \imath }{\sqrt{2}} \int_{-\infty}^{+\infty} \psi^\ast
 \dv{\psi}{x} dx,
$ (real quantity  via an integration by part).

It is very convenient to introduced the annihilation operator   $\ba$  and creation operator  $\ba^\dag$:
\begin{equation*}\label{eq:aa}
 \ba=\bX+ \imath  \bP\equiv\frac{1}{\sqrt{2}}
 \left(x+\dv{}{x}\right)
 ,\quad
 \ba^\dag=\bX- \imath \bP \equiv \tfrac{1}{\sqrt{2}}
 \left(x-\dv{}{x}\right)
 .
\end{equation*}
We have
$$
 [\bX,\bP] = \tfrac{ \imath  }{2}\bI, \quad  [\ba,\ba^\dag]= \bI, \quad
  \bH = \omega (\bP^2 + \bX^2)  =\omega \left(\ba^\dag \ba + \tfrac{1}{2} \bI \right)
$$
where $\bI$ stands for the identity operator.

Since   $[{\ba},{\ba^\dag}]= \bI$,  the spectral decomposition of ${\ba^\dag} {\ba} $
 is simple.  The Hermitian operator $\bN=\ba^\dag\ba$, the photon-number operator,   admits  $\NN$  as non degenerate spectrum. The normalized eigenstate $\ket{n}$ associated to $n\in\NN$,  is denoted by  $\ket{{n}}$. Thus the underlying  Hilbert  space reads
 $$\HH=\left\{\sum_{n\ge 0} \psi_n \ket{n},\; (\psi_n)_{n\ge 0}\in l^2(\CC)\right\}$$ where $(\ket{n})_{n\in\NN}$ is the  Hilbert basis of photon-number states (also called Fock states).
 For $n>0$, we have
$$
 {\ba}{\ket{{n}}} = \sqrt{n} ~{\ket{{n-1}}}
 , \quad
 {\ba^\dag} {\ket{{n}}} = \sqrt{n+1}~{\ket{{n+1}}}
 .
$$
The ground state  ${\ket{{0}}}$  is characterized by    ${\ba}{\ket{{0}}}= 0$.  It corresponds to the Gaussian function
${\psi_0(x)}=\frac{1}{\pi^{1/4}}\exp(-x^2/2)$.

For any function  $f$  we have the following commutations
$$
\ba f(\bN)= f(\bN+\bI) \ba, \quad \ba^\dag  f(\bN)= f(\bN-\bI) \ba^\dag.
$$
In particular for any angle $\theta$,  $e^{i\theta \bN} \ba e^{-i\theta \bN} = e^{-i\theta} \ba $.

For any amplitude $\alpha\in\CC$,  the Glauber displacement unitary operator $\bD_\alpha$  is defined by
\begin{equation*}\label{eq:Dalpha}
    \bD_\alpha = e^{\alpha~ \ba^\dag - \alpha^*\ba}
\end{equation*}
We have $\bD_\alpha^{-1}= \bD_\alpha^\dag = \bD_{-\alpha}$.
The following Glauber formula is useful: if two operators $\bA$ and  $\bB$ commute with their commutator,    i.e., if   $[\bA,[\bA,\bB]]=[\bB,[\bA,\bB]]=0$, then we have
$e^{\bA+\bB} =e^{\bA} ~e^{\bB}~ e^{-\tfrac{1}{2}[\bA,\bB]}$.
Since  $\bA=\alpha \ba^\dag$ and   $\bB=-\alpha^* \ba$ are in this case, we have another expression for  $\bD_\alpha$
\begin{equation*}\label{eq:DalphaBis}
    \bD_\alpha = e^{-\tfrac{|\alpha|^2}{2}}~e^{\alpha \ba^\dag} e^{ - \alpha^* \ba}= e^{+\tfrac{|\alpha|^2}{2}}~e^{-\alpha^* \ba} e^{  \alpha \ba^\dag}
    .
\end{equation*}
The terminology displacement has its origin in the following property derived from Baker-Campbell-Hausdorff formula:
\begin{equation*}\label{eq:DalphaA}
\forall \alpha\in\CC, \quad \bD_{-\alpha} \ba \bD_\alpha = \ba + \alpha \quad  \text{and} \quad \bD_{-\alpha} \ba^\dag \bD_\alpha = \ba^\dag + \alpha^*
.
\end{equation*}
To the classical state $(x,p)$ is associated a  quantum  state usually called coherent state of complex amplitude $\alpha= (x+\imath p)/\sqrt{2}$ and denoted by   $\ket{\alpha}$:
\begin{equation}\label{eq:coherent}
    \ket{\alpha} = \bD_\alpha \ket{0} = e^{-\tfrac{|\alpha|^2}{2}} \sum_{n=0}^{+\infty} \tfrac{\alpha^n}{\sqrt{n!}} \ket{n}
    .
\end{equation}
$\ket{\alpha}$  corresponds to the translation of the Gaussian profile  corresponding to vacuum state  $\ket{0}$:
$$
 \ket{\alpha} \equiv \left( \RR\ni x \mapsto  \tfrac{1 }{\pi^{1/4}}
e^{\imath\sqrt{2}x \Im\alpha } e^{-\frac{(x-\sqrt{2}\Re\alpha )^2}{2}}\right).
$$
This usual notation is potentially ambiguous: the coherent state $\ket{\alpha}$ is very different from  the photon-number state $\ket{n}$ where $n$ is a non negative  integer:
The probability  $p_n$  to obtain  $n\in\NN$  during the measurement of  $\bN$ with  $\ket\alpha$ obeys to a Poisson law  $p_n=e^{-|\alpha|^2} |\alpha|^{2n}/n!$. The resulting  average energy   is thus given by $\bket{\alpha|\bN|\alpha}= |\alpha|^2$.
Only for $\alpha=0$ and $n=0$, these quantum states coincide.

The coherent state $\alpha\in\CC$  is the unitary eigenstate of  $\ba$ associated to the eigenvalue  $\alpha\in\CC$:
$\ba \ket{\alpha} =\alpha \ket{\alpha}$.
Since  $\bH/\hbar =\omega (\bN+\tfrac{1}{2})$, the solution of the Schr\"{o}dinger equation
$\dotex \ket\psi = - \imath \frac{\bH}{\hbar } \ket\psi,$ with initial value  a coherent state $\ket{\psi}_{t=0}=\ket{\alpha_0}$ ($\alpha_0\in\CC$)   remains a coherent state with time varying  amplitude $\alpha_t= e^{-\imath \omega t}\alpha_0$:
$$
\ket\psi_t= e^{-\imath \omega t/2} \ket{\alpha_t}
.
$$
These coherent  solutions  are the quantum counterpart of the classical  solutions: $x_t=\sqrt{2}\Re(\alpha_t)$ and $p_t=\sqrt{2}\Im(\alpha_t)$ are solutions of the classical Hamilton equations $\dotex x = \omega p$ and $\dotex p = -\omega x$ since  $\dotex \alpha_t = - \imath \omega \alpha_t $.
The addition of a control input,  a classical  drive  of amplitude  $u\in\RR$,   yields to the following control Schr\"{o}dinger equation
$$
\dotex \ket\psi = -\imath \Big(\omega \left(\ba^\dag \ba + \tfrac{1}{2}\right) + u (\ba+\ba^\dag) \Big) \ket\psi
$$
It is the quantum version of the control classical harmonic oscillator
$$
\dotex x = \omega p, \quad \dotex p = -\omega x -  u \sqrt{2}
.
$$
\begin{figure}
  \centerline{ \includegraphics[width=0.8\textwidth]{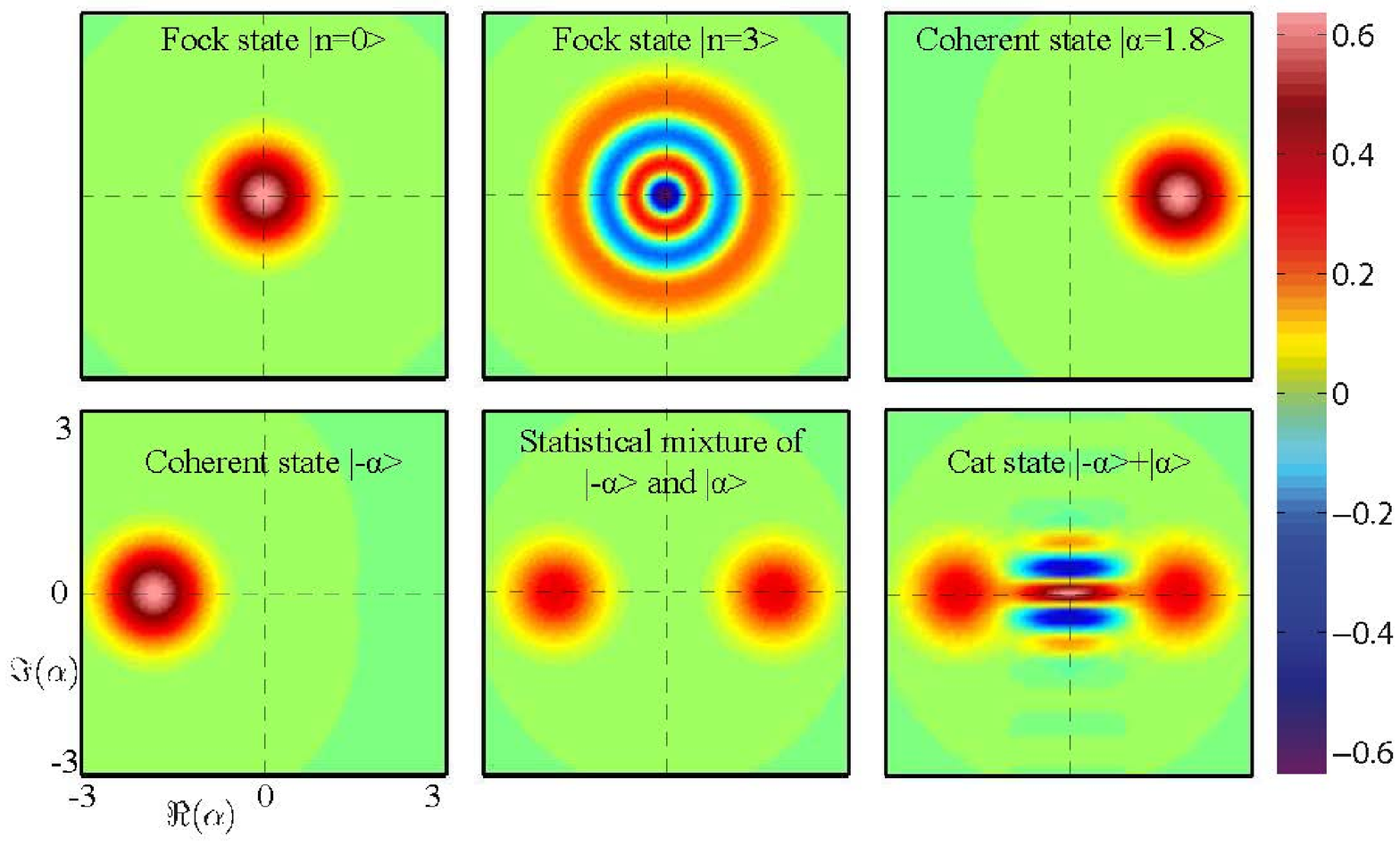}}
  \caption{Wigner function of typical quantum states of an harmonic oscillator. }\label{fig:Wigner_intro}
\end{figure}
A possible definition of the Wigner function $W^\rho$  attached to any density operator $\rho$  is as follows:
$$
  W^\rho:\CC\ni\alpha\rightarrow \tfrac{2}{\pi}\tr {e^{i\pi\bN} e^{-\alpha \ba^\dag + \alpha^* \ba} \rho  e^{\alpha \ba^\dag - \alpha^* \ba} } \in[-2/\pi,2/\pi]
$$
where $\alpha=x+i p$ is a  position in the phase-plane $(x,p)$ of the classical oscillator. With the correspondences
\begin{equation}
\begin{split}\label{eq:correspondence}
&  \dv{}{\alpha}= \tfrac{1}{2} \left(\dv{}{x}-i \dv{}{p} \right), \quad \dv{}{\alpha^*}= \tfrac{1}{2} \left(\dv{}{x}+i \dv{}{p} \right)
\\
&W^{\rho \ba } = \left(\alpha - \tfrac{1}{2} \dv{}{\alpha^*}\right) W^\rho, \quad  W^{\ba \rho } = \left(\alpha + \tfrac{1}{2} \dv{}{\alpha^*}\right) W^\rho
\\
&W^{\rho \ba^\dag } = \left(\alpha^* + \tfrac{1}{2} \dv{}{\alpha}\right) W^\rho, \quad  W^{\ba^\dag \rho } = \left(\alpha^* - \tfrac{1}{2} \dv{}{\alpha}\right) W^\rho
\end{split}
\end{equation}
the Lindblad-Kossakovki governing the evolution of the density operator $\rho$  of  a  quantum oscillator,  with damping time constant $1/\kappa >0$ and  resonant  drive of real amplitude $u$,
$$
    \dotex \rho = u[\ba^\dag - \ba ,\rho] + \kappa\left(\ba\rho \ba^\dag - (\bN \rho + \rho \bN)/2\right)
    ,
$$
becomes a convection-diffusion equation for the Wigner function  $W^\rho$
$$
 \dv{W^\rho}{t} =
    \tfrac{\kappa}{2} \left( \dv{~}{x} \Big( (x-\bar\alpha) W^\rho\Big) +  \dv{~}{p} \Big( p W^\rho\Big)
    + \tfrac{1}{4}  \Delta W^\rho \right)
$$
where  $\Delta$ denotes the Laplacian operator $\tfrac{\partial ^2}{\partial x^2}+ \tfrac{\partial ^2}{\partial p^2}$.
The  solutions converge  toward the Gaussian steady-state   $W^{\rhob}(x,p) = \tfrac{2}{\pi} e^{-2 (x-\bar\alpha)^2-2p^2}$, where $\rhob=\ket{\bar\alpha}\bra{\bar\alpha}$  is the coherent  state of amplitude $\bar\alpha= 2u/\kappa$.

\section{Qubit} \label{ap:qubit}

The underlying Hilbert space $
\HH=\CC^2=\left\{c_g\ket{{g}}+c_e\ket{{e}},\; c_g,c_e\in \CC\right\}$
where $(\ket{g},\ket{e})$ is the  ortho-normal frame formed by the ground state $\ket{g}$ and the excited state $\ket{e}$.
It is usual to consider the following operators on $\mathcal H$:
\begin{equation} \label{eq:Smpxyz}
  \begin{split}
    &\bSm=\ket{ {g}}\bra{ {e}}, \quad \bSp=\bSm^\dag=\ket{ {e}}\bra{ {g}}, \quad {\bSx}=\bSm+\bSp=\ket{ {g}}\bra{ {e}}+ \ket{ {e}}\bra{ {g}},
    \\
 & \bSy=i \bSm-i \bSp = i\ket{ {g}}\bra{ {e}}-i\ket{ {e}}\bra{ {g}}
 ,\quad  {\bSz}=\bSp\bSm - \bSm \bSp = \proj{e}-\proj{g}.
  \end{split}
\end{equation}
$\bSx$, $\bSy$ and $\bSz$ are the Pauli operators. They are square root of $\bI$:
$
\bSx^2=\bSy^2= \bSz^2=\bI.
$
They anti-commute
$$
\bSx \bSy= - \bSy \bSx = i \bSz, \quad
\bSy \bSz= - \bSz \bSy = i \bSx, \quad
\bSz \bSx= - \bSx \bSz = i \bSy
$$
and thus $ [\bSx,\bSy] = 2 i \bSz$,  $[\bSy,\bSz] = 2 i \bSx$,  $[\bSz,\bSx] = 2 i \bSy$.
The uncontrolled  evolution is governed by the  Hamiltonian $\bH/\hbar={\omega \bSz}/2$ where $\omega >0$ is the qubit pulsation. Thus the solution of $\dotex \ket\psi = -i\frac{\bH}{\hbar}\ket\psi$ is given by
$$
\ket\psi_t = e^{-i\left(\tfrac{\omega t}{2}\right)\bSz} \ket\psi_0= \cos\left(\tfrac{\omega t}{2}\right) \ket\psi_0- i\sin\left(\tfrac{\omega t}{2}\right) \bSz\ket\psi_0
$$
since for any angle $\theta$ we have
$$
e^{i\theta\bSx}= \cos\theta + i \sin\theta\bSx, \quad
 e^{i\theta\bSy}= \cos\theta + i \sin\theta\bSy,\quad
 e^{i\theta\bSz}= \cos\theta + i \sin\theta\bSz.
$$
Since the Pauli operators anti-commute,  we have  the useful  relationships:
$$
 e^{i\theta\bSx}\bSy =  \bSy e^{-i\theta\bSx}, \quad
 e^{i\theta\bSy}\bSz =  \bSz e^{-i\theta\bSy}, \quad
 e^{i\theta\bSz}\bSx =  \bSx e^{-i\theta\bSz}.
$$

The orthogonal projector  $\rho =\ket\psi \bra\psi$, the density operator associated to the pure state $\ket\psi$,  obeys to the Liouville  equation $\dotex \rho = -\tfrac{i}{\hbar} [\bH,\rho].$ Mixed quantum states are described by  $\rho$ that are Hermitian, non-negative and of  trace one. For a qubit, the Bloch sphere representation is a useful tool exploiting the  smooth correspondence between such $\rho$ and the unit ball of $\RR^3$  considered as Euclidian space:
$$
\rho = \frac{\bI + x\bSx + y\bSy + z \bSz}{2}, \quad (x,y,z) \in\RR^3, \quad x^2+y^2+z^2\leq 1
.
$$
 $(x,y,z)\in\RR^3$  are  the coordinates in the orthonormal frame $(\vec\imath,\vec\jmath,\vec k)$   of the Bloch vector $\vec M\in\RR^3$. This vector  lives on or inside the unit sphere, called Bloch sphere:
$$
\vec M = x \vec\imath + y \vec\jmath +  z \vec k.
$$
Since $\tr{\rho^2}=x^2+y^2+z^2$, $\vec M$ is on the Bloch sphere when $\rho$ is of rank one and thus is a pure state.
The translation of Liouville equation on $\vec M$ yields with $\bH/\hbar = \omega \bSz/2$:
$
\dotex \vec M = \omega \vec k \times \vec M.
$
For the two-level system with the coherent  drive  described by the complex-value control $u$,  $\bH/\hbar = \tfrac{\omega}{2} \bSz + \tfrac{\Re(u)}{2} \bSx + \tfrac{\Im(u)}{2} \bSy$ and the Liouville equation reads, with the Bloch vector $\vec M$ representation,
$$
\dotex \vec M = (\Re(u)\vec\imath +\Im(u) \vec\jmath+\omega \vec k ) \times \vec M
.
$$

\section{Jaynes-Cumming Hamiltonians and propagators} \label{ap:JC}

The Jaynes-Cummings Hamiltonian~\cite{jaynes-cummings-ieee64}  is the simplest Hamiltonian describing  the interaction between an harmonic oscillator and a qubit.  Such an interaction admits two regimes, the resonant one where the oscillator and the qubit  exchange energy, the dispersive one where the oscillator pulsation depends on the qubit state and where the qubit pulsation, slightly different from the oscillator pulsation, depends on the number of vibration quanta.  We recall below the simplest forms of these   Hamiltonians in the  interaction frame. A  deeper and complete  presentation   can be found in~\cite{haroche-raimondBook06}.

The resonant Hamiltonian $\bH_{res}$ is given by
\begin{equation}\label{eq:Hres}
  \bH_{res} /\hbar = i f(t)  ~\big( \ba^\dag\otimes\bSm - \ba\otimes \bSp \big)=  i f(t)  ~\big( \ba^\dag\otimes\ket{g}\bra{e} - \ba\otimes \ket{e}\bra{g} \big)
\end{equation}
whereas the dispersive one $\bH_{disp}$ is a simple tensor product:
\begin{equation}\label{eq:Hdisp}
  \bH_{disp} /\hbar =f(t)~\bN \otimes \bSz= f(t) ~\bN\otimes \big( \proj{e}-\proj{g}\big)
\end{equation}
where $f(t)$ is a  known real  parameter depending possibly   on the  time  $t$.

Simple computations show that the resonant propagator $\bU_{res}$ between $t_0$ and $t_1$  associated to $\bH_{res}$, i.e., the solution of  Cauchy problem
$$
\dotex \bU = - i \frac{\bH_{res}}{\hbar} \bU, \quad U(t_0)=\bI,
$$
is  explicit and given by the following compact formulae:
\begin{multline}\label{eq:Ures}
  \bU_{res}(t_0,t_1) = \cos\left( \tfrac{\int_{t_0}^{t_1}f}{2} \sqrt{\bN}\right) \otimes\ket{g}\bra{g} +  \cos\left( \tfrac{\int_{t_0}^{t_1}f}{2}\sqrt{\bNp}\right)\otimes \ket{e}\bra{e} \\
 - \ba \frac{\sin\left( \tfrac{\int_{t_0}^{t_1}f}{2} \sqrt{\bN}\right)}{\sqrt{\bN}} \otimes \ket{e}\bra{g} \, +  \frac{\sin \left( \tfrac{\int_{t_0}^{t_1}f}{2} \sqrt{\bN}\right)}{\sqrt{\bN}} \, \ba^\dag \otimes\ket{g}\bra{e}
 .
\end{multline}
It is instructive to check  that $ \bU^\dag_{res}\bU_{res}=\bI$.
Similarly, the dispersive  propagator  $\bU_{disp}$  between $t_1$ and $t_2$  associated to $\bH_{disp}$ is given by
\begin{equation}\label{eq:Udisp}
  \bU_{disp}(t_0,t_1) = \exp\left(i \bN \int_{t_0}^{t_1}\hspace{-0.7em}f \right)\otimes \ket{g}\bra{g} +\exp\left(-i \bN \int_{t_0}^{t_1}\hspace{-0.7em}f \right) \otimes \ket{e}\bra{e}.
\end{equation}

\section{A positiveness-preserving numerical scheme} \label{ap:num}
This appendix describes  a positiveness-preserving formulation  of the Euler-Milstein  scheme for the numerical integration of    stochastic master equations driven by a single Wiener process. They admit  the following  form
\begin{multline} \label{eq:lindblad}
    d\rho_t = \left(-\imath[\bH,\rho_t] + \sum_{\mu} \bL_\mu \rho_t \bL_\mu^\dag  - \tfrac{1}{2} (\bL_\mu^\dag \bL_\mu \rho_t + \rho_t \bL_\mu^\dag \bL_\mu)  \right) dt
    \\
    + \bigg(
     \bL \rho_t \bL^\dag  - \tfrac{1}{2} (\bL^\dag \bL \rho_t + \rho_t \bL^\dag \bL)\bigg) dt +
     \sqrt{\eta} \bigg( \bL\rho_t + \rho_t \bL^\dag - \tr{ \bL\rho_t + \rho_t \bL^\dag}\rho_t \bigg) dW_t
\end{multline}
where $\rho$ is a square  non-negative Hermitian matrix of trace $1$, $\bL_\mu$ and  $\bL$ are  square  matrices, $ W_t$ is a Wiener process and $\eta\in[0,1]$ is the detection efficiency. The measured continuous signal $y_t$  is given by $dy_t =   \sqrt{\eta} \tr{\bL\rho_t + \rho_t \bL^\dag} dt + dW_t $.

For $dx =f(x)dt + g(x) dW_t$ ($x\in\RR^d$ for some integer $d$, $f$ and $g$ smooth functions),  the Euler-Milstein scheme (order $1$ in  the discretization step denoted $dt$) reads \cite{MilsteinBook}
$$
x_{n+1}= x_n + f(x_n) dt  + g(x_n) dW_n + \tfrac{1}{2}  \frac{\partial g}{\partial x}(x_n) \cdot g(x_n) (dW_n^2-dt)
$$
where $x_{n}$, for $n\in \NN$,  is an approximation of $x_{n dt}$ and $dW_n$ is a Gaussian variable with zero average and  variance $dt$.
For~\eqref{eq:lindblad}, we get
\begin{multline*}
  \rho_{n+1} = \rho_n +
  \left(-\imath[\bH,\rho_n] + \sum_{\mu} \bL_\mu \rho_n \bL_\mu^\dag  - \tfrac{1}{2} (\bL_\mu^\dag \bL_\mu \rho_n + \rho_n \bL_\mu^\dag \bL_\mu)  \right) dt
    \\+\bigg( \bL \rho_n \bL^\dag  - \tfrac{1}{2} (\bL^\dag \bL \rho_n + \rho_n \bL^\dag \bL)\bigg) dt + \sqrt{\eta} \bigg( \bL\rho_n + \rho_n \bL^\dag
    - \tr{ \bL\rho_n + \rho_n\bL^\dag}\rho_n \bigg) dW_n \\
    +
    \tfrac{\eta}{2} \bigg(
      \bL^2 \rho_n + \rho_n (\bL^\dag)^2 + 2 \bL\rho_n \bL^\dag - 2 \tr{\bL\rho_n+\rho_n \bL^\dag} (\bL\rho_n+\rho_n \bL^\dag)  \ldots
      \\
      - \big( \tr{ \bL^2 \rho_n + \rho_n (\bL^\dag)^2} + 2 \tr{\bL\rho_n \bL^\dag} - 2 \trr{\bL\rho_n + \rho_n \bL^\dag}  \big)\rho_n
       \bigg) (dW_n^2-dt)
       .
\end{multline*}
Let us consider the following matrix
$$
\bM_n = \bI - dt \left( i \bH +\tfrac{1 }{2}\sum_{\mu}  \bL_\mu^\dag \bL_\mu +  \tfrac{1 }{2}\bL^\dag \bL \right)   + \sqrt{\eta} \bigg(\sqrt{\eta}  \tr{\bL\rho_n +\rho_n\bL^\dag} dt + dW_n \bigg) \bL
+ \tfrac{\eta}{2} (dW_n^2-dt)  \bL^2.
$$
Here  $dW_n$ is  of order $\sqrt{dt}$ and  $dW_n^2-dt$ is  of order $dt$. Then simple but slightly tedious computations up to $dt^{3/2} $ show that $\rho_{n+1}$ given by the above  Euler-Milstein scheme   reads also
\begin{equation}\label{eq:Milstein}
\rho_{n+1} = \frac{\bM_n \rho_n \bM_n^\dag + \sum_\mu  \bL_\mu \rho_n \bL_\mu^\dag dt  +  (1-\eta) \bL\rho_n \bL^\dag dt }{\tr{
\bM_n \rho_n \bM_n^\dag +\sum_\mu  \bL_\mu \rho_n \bL_\mu^\dag dt  +   (1-\eta)   \bL\rho_n \bL^\dag dt} } + O(dt^{3/2}).
\end{equation}
 When  $\eta=0$, this expression   provides, for any deterministic Lindblad differential equation,    a positiveness-preserving formulation of the  explicit Euler scheme.

\end{document}